\documentclass[12pt]{article}
\usepackage{amssymb,latexsym,
amsfonts,amsbsy}

\newcommand{\qed}{\hbox{\rule[-2pt]
{3pt}{6pt}}}

\newtheorem{dfe}{Definition}
[section]

\newtheorem{theo}[dfe]{Theorem}

\newtheorem{pro}[dfe]{Proposition}

\makeatletter

\@addtoreset{equation}{section}
\makeatother



\begin{title}
{Euclidean designs and coherent configurations}
\end{title}
\author{Eiichi Bannai\quad and \quad
Etsuko Bannai\\
}

\begin{document}
\maketitle

\begin{abstract}
The concept of spherical $t$-design, which is a finite subset of the 
unit sphere, was introduced by Delsarte-Goethals-Seidel (1977). 
The concept of Euclidean $t$-design, which is a two step generalization of spherical design in the sense that it is a finite weighted subset of Euclidean space, by Neumaier-Seidel (1988). 
We first review these two concepts, as well as the concept of tight $t$-design, i.e., the one whose cardinality reaches the 
natural lower bound.  We are interested in $t$-designs (spherical or 
Euclidean) which are either tight or close to tight.  
As is well known by
Delsarte-Goethals-Seidel (1977), in the study of spherical $t$-designs and in particular of those which are either tight or close to tight, association schemes play important roles. The main 
purpose of this paper is to show that in the study of Euclidean $t$-designs and in particular of those which are either tight or close to tight, coherent configurations play important roles. Here, 
coherent configuration is a purely combinatorial concept
defined by D. G. Higman, and is obtained by axiomatizing 
the properties of general, not necessarily transitive,
permutation groups, in the same way as association scheme was obtained by axiomatizing the properties of transitive permutation groups. The main purpose of this paper is to prove that Euclidean 
$t$-designs satisfying certain conditions give the structure of 
coherent configurations. In particular, it is seen that a tight 
Euclidean $t$-design on two concentric spheres centered at the 
origin has the structure of coherent configuration. Moreover, 
as an
application of this general theory, we discuss the current status
of our research to try to classify Euclidean $4$-designs $(X,w)$
on two concentric spheres $S=S_1\cup S_2$ centered at the origin
whose weight function is constant on
each $X\cap S_i\ (i=1,2)$ and the number of the
inner products between the distinct two points in 
$X\cap S_i$ and $X\cap S_j$ is at
most $2$ for $i, j=1,2$. We describe all the parameters of the 
coherent configurations, in terms of the parameters of the Euclidean 
designs. 
The classification of such Euclidean $4$-designs is not yet completed, but we have found two new families of feasible
parameters of such Euclidean $4$-designs and the associated  coherent configurations. One family corresponds to 
Euclidean tight $4$-designs on two concentric spheres 
and another family is obtained from non-tight Euclidean 
4-designs (and is related to the spherical tight $4$-designs of 
one dimension more). 
\end{abstract} 

\section{Introduction}
Spherical $t$-designs are defined in the paper
by Delsarte-Goethals-Seidel \cite{D-G-S}.
In that paper they showed that spherical designs satisfying some conditions have structures of 
Q-polynomial association scheme.
Euclidean $t$-designs are defined in the paper by
Neumaier-Seidel \cite{N-S} as a generalization
of spherical designs. There are very natural lower bounds
for the cardinalities of Euclidean $t$-designs
 (see \cite{M,D-S,N-S,B-1,B-B-5}) and tightness for the Euclidean
designs defined (see \cite{D-S,N-S,B-B-3,B-B-5,B-1,B-B-H-S}). It is an interesting and important problem
to construct and classify Euclidean tight designs.
Examples of tight Euclidean $t$-designs 
are constructed (see \cite{Baj-1,Baj-2, B-B-1,
B-B-3,B-B-H-S,B-1,B-2}).
We observed that some of the examples of tight 
Euclidean $t$-designs constructed have the structures of
coherent configurations. Coherent configuration is a
concept defined by Higman \cite{H-1, H-2} as a generalization of
association schemes. In this paper
we give sufficient conditions for Euclidean designs
to have the structures of coherent configurations.
In particular we prove that the Euclidean tight $t$-designs
supported by two concentric spheres have the
structures of coherent configurations.
We give a series of feasible parameters for 
Euclidean $4$-designs supported by 2 concentric spheres
having the structures of coherent configurations
and we also give a series of feasible parameters for
tight Euclidean $4$-designs supported by 2
concentric spheres. 

First we give some notation. 
Let $X$ be a finite set in Euclidean space $\mathbb R^n$. 
Let $w$ be a positive real valued weight function
defined on $X$. 
We assume
$n\geq 2$ throughout this paper and consider the
weighted finite sets $(X,w)$ in $\mathbb R^n$. 
Let $\boldsymbol x\cdot\boldsymbol y$
be the canonical inner product between
$\boldsymbol x$ and $\boldsymbol y$ in $\mathbb R^n$ and $\|\boldsymbol x\|=\sqrt{\boldsymbol x
\cdot \boldsymbol x}$.
Let $S^{n-1}$ be the unit sphere centered at the origin.
Let $S^{n-1}(r)$ be the sphere of 
radius $r$ centered at the origin, where
$r$ possibly be $0$.
We can decompose $X$ into a disjoint union 
of nonempty subsets in the following manner,
that is, $r_1,\ r_2,\ldots ,r_p$ are 
distinct nonnegative
real numbers and
$X=X_1\cup X_2\cup\cdots \cup X_p,$
$X_i\subset S^{n-1}(r_i)$ for $i=1,\ 2,\ldots, p$. Let us denote $S_i= S^{n-1}(r_i),\ 1\leq i\leq p$.
Let $S=\cup_{i=1}^pS_i$. 
Let $\varepsilon_S=1$ if $0\in X$, i.e., if there exists
$i$ satisfying $r_i=0$ and $\varepsilon_S=0$, otherwise. 
We say
$X$ is supported by $p$ concentric spheres.
Let $w(X_i)=\sum_{\boldsymbol x\in X_i}
w(\boldsymbol x)$ for $i=1,\ldots, p$.
Let $\sigma$ and $\sigma_i,$ $1\leq i\leq p$, be the 
Haar measure on $S^{n-1}$ and $S_i,$ $1\leq i\leq p$,
respectively. Let $|S^{n-1}|=\int_{S^{n-1}}d\sigma
(\boldsymbol x)$, $|S_i|=\int_{S_i}d\sigma_i
(\boldsymbol x)$, $1\leq i\leq p$.
Here, if $r_i=0$, then we define
$\frac{1}{|S_i|}\int_{S_i}f(\boldsymbol x)d\sigma_i
(\boldsymbol x)=f(0)$ for any polynomial
 $f(\boldsymbol x)$.
We assume $|S_i|={r_i}^{n-1}|S^{n-1}|$ for $r_i>0$. 
$\mathcal P(\mathbb R^n)$ denotes the vector 
space of polynomials
in $n$ variables $x_1,\ldots, x_n$ over the 
fields $\mathbb R$ of real numbers.
Let
$\mbox{Hom}_l(\mathbb R^n)$ be the
subspace of  ${\cal P}(\mathbb R^n)$
which consists of homogeneous
polynomials of degree $l$.
Let ${\cal P}_l(\mathbb R^n)
=\oplus_{i=0}^l
\mbox{Hom}_i(\mathbb R^n)$.
Let $\mbox{Harm}(\mathbb R^n)$ be
the subspace of
${\cal P}(\mathbb R^n)$
which consists of all the
harmonic polynomials. Let
$\mbox{Harm}_l(\mathbb R^n)
=\mbox{Harm}(\mathbb R^n)\cap
\mbox{Hom}_l(\mathbb R^n).$\\

The following is the definition of Euclidean $t$-design.
See Remarks after Theorem \ref{theo:N-S} also.
\begin{dfe}[Euclidean $\boldsymbol t$-design]{\rm(see \cite{N-S})}
\label{dfe:1-1}\quad
Let $t$ be a natural
number. A weighted finite set $(X,w)$ in $\mathbb R^n$ is a 
Euclidean $t$-design,
if the following
equation
$$\sum_{i=1}^p\frac{w(X_i)}{|S_i|}
\int_{\boldsymbol x\in
S_i}f(\boldsymbol x)d\sigma_i(\boldsymbol x)
=\sum_{\boldsymbol u\in X}w(\boldsymbol u)
f(\boldsymbol u)$$
is satisfied for any polynomial 
$f\in \mathcal P_t(\mathbb R^n)$.
\end{dfe}
{\bf Remark:} 
If $r>0$, $X\subset S^{n-1}(r)$ and
$\frac{1}{r}X(\subset S^{n-1})$ is a spherical  
$t$-design, then we also call $X$ a spherical $t$-design. With this definition, if $p=1$, $X\not=\{0\}$,
and $w(\boldsymbol x)\equiv 1$ in Definition
\ref{dfe:1-1}, then
$X$ is a spherical $t$-design.\\

For the cardinalities of Euclidean $t$-designs, natural lower bounds are proved by 
M\"oller in 1978 (see \cite{M-0,M}, also \cite{B-1,B-B-5,D-S,N-S})
and concept of tightness are defined.
Here we give only for the case where $t$ is even.
The definition of tightness for $t$ odd is more delicate.
Reader can find more detailed information 
in \cite{B-B-5,B-B-H-S,M}. 

\begin{theo}[\cite{M,D-S}]\label{theo:lowerbd}
Let $(X,w)$ be a Euclidean $2e$-design supported by $p$ concentric spheres $S$ in $\mathbb R^n$.
Then 
$$|X|\geq \dim(\mathcal P_e(S))$$
holds, where $\mathcal P_e(S)=\{f|_S\mid 
f\in \mathcal P_e(\mathbb R^n)\}$.
\end{theo}

\begin{dfe}[\cite{D-S,B-B-3}]\label{dfe:tight}
\begin{enumerate}
\item Definition and notation are the same as above.
If equality holds in Theorem \ref{theo:lowerbd},
then $(X,w)$ is called a tight $2e$-design on $p$ concentric spheres
\item Moreover if $\dim(\mathcal P_e(S))
=\dim(\mathcal P_e(\mathbb R^n))
(={n+e\choose e})$ holds, then
$(X,w)$ is called a Euclidean tight $2e$-design.
\end{enumerate}
\end{dfe}

We give some more notation. Let $(X,w)$ be a Euclidean $t$-designs supported by
$p$ concentric spheres.
For any $X_\lambda,X_\mu\not=\{0\}$, let 
$$A(X_\lambda,X_\mu)=A(X_\mu,X_\lambda)=\left\{\frac{\boldsymbol x\cdot\boldsymbol y}
{r_\lambda r_\mu}\ \bigg|\ \boldsymbol x\in X_\lambda,
\boldsymbol y\in X_\mu,
 \boldsymbol x\not= \boldsymbol y \right\}.$$
Let $s_{\lambda,\mu}=s_{\mu,\lambda}=|A(X_\lambda,X_\mu)|$
and $A(X_\lambda,X_\mu)=\{\alpha_{\lambda,\mu}^{(u)}=
\alpha_{\mu,\lambda}^{(u)}
\mid u=1,\ldots,s_{\lambda,\mu}\}$.
Let $\alpha_{\lambda,\lambda}^{(0)}=1$ for any 
$X_\lambda
\not=\{0\}$.\\

The following are the main theorems of this paper.
\begin{theo}\label{theo:1-2}
Let $(X,w)$ be a Euclidean $t$-design in 
$\mathbb R^n$ supported by $p$ concentric spheres.
Assume $w(\boldsymbol x)\equiv w_\nu$
for any $\boldsymbol x\in X_\nu$ $(1\leq \nu\leq p)$.
Moreover we assume the following {\rm(1)} 
or {\rm(2)}.
\begin{enumerate}
\item If
$s_{\lambda,\nu}+s_{\nu,\mu}\leq t-2(p-\varepsilon_S-2)$ holds
for any $\lambda$, $\nu$ and $\mu$
with $1\leq \lambda,\nu,\mu\leq p$.
\item 
If $X$ is antipodal and $s_{\lambda,\nu}+s_{\nu,\mu}-\delta_{\lambda,\nu}
-\delta_{\nu,\mu}\leq t-2(p-\varepsilon_S-2)$ holds for any
$\lambda,\ \nu$ and $\mu$ satisfying
$1\leq \lambda,\nu,\mu\leq p$.
\end{enumerate}
Then $X$ has the structure of a coherent 
configuration. 
\end{theo}
\begin{theo}\label{theo:1-3}
Let $t\geq 2$ and
$(X,w)$ be a tight Euclidean $t$-design
supported by 2 concentric spheres.
Then $X$ has the structure of a coherent configuration.
\end{theo}
{\bf Remark:} If $t=1$, then $X$ consists of an antipodal pair in $\mathbb R^n$ and $p=1$. 
\begin{theo}\label{theo:1-4}
Let $(X,w)$ be a Euclidean $4$-design in
$\mathbb R^n$
supported by 2 concentric spheres. Assume
$0\not\in X$,
$w$ is constant on each $X_\lambda$,
and $s_{\lambda,\mu}\leq 2$  ($\lambda,\mu=1,2$).
Then $X$ has the structure of a coherent configuration
and the following holds.
\begin{enumerate}
\item $s_{1,2}=2$.
\item $(X,w)$ is either a tight Euclidean $4$-design or similar to one of the Euclidean $4$-designs having the following parameters. 

{\rm(i)} $n=2$, 
$X_1=\{\pm(\frac{1}{\sqrt{2}},\ \frac{1}{\sqrt{2}}),\ 
\pm(\frac{1}{\sqrt{2}},\ -\frac{1}{\sqrt{2}})\}$,
$X_2=\{(\pm r_2,\ 0),\ (0,\ \pm r_2)\}$, 

$w(\boldsymbol x)=1$, for any $\boldsymbol x\in X_1$
and $w(\boldsymbol x)=r_2^{-4}$ for $\boldsymbol x\in X_2$,
where $r_2$ is any positive real number $r_2\not=1$.\\

{\rm(ii)} $n=(2k-1)^2-4$, where $k$ is any integer satisfying $k\geq 2$,

$|X_1|=2(2k+1)(k-1)^3$,   
$|X_2|=2k^3(2k-3)$, 

$A(X_1)=\{\frac{k-2}{k(2k-3)},\
-\frac{1}{2k-3}\}$, $A(X_2)=\{\frac{1}{2k+1},\
-\frac{k+1}{(k-1)(2k+1)}\}$,
$A(X_1,X_2)=\{\frac{1}{\sqrt{n}},\
-\frac{1}{\sqrt{n}}\}$, 
$r_1=1$, $w(\boldsymbol x)=1$ for $\boldsymbol x
\in X_1$ and
$w(\boldsymbol x)= \frac{(2k+1)^2(k-1)^4}{(2k-3)^2k^4}r_2^{-4}$, where $r_2$ is any positive real number satisfying $r_2\neq 1$.

The intersection numbers of the
corresponding coherent configurations are given 
as polynomials of $k$ (see 
Appendix I).\end{enumerate}
\end{theo}
{\bf Remark:}\\
(1) In Theorem \ref{theo:1-4}, if $|X_1|=n+1$, then
$X_1$ must be a tight spherical $2$-design, i.e.,
a regular simplex on $S_1$. Also we will prove,
in \S 4.1 (Theorem \ref{theo:4-1}), that in this case 
$(X,w)$ must be a tight Euclidean 4-design. 
Tight Euclidean 4-designs
with this property are classified in \cite{B-2}.\\
(2) Let $r_2=\frac{k-1}{k}
\sqrt{\frac{2k+1}{2k-3}}$ in 
Theorem \ref{theo:1-4} (2)(ii), then the corresponding Euclidean 4-design
$(X,w)$ is of constant weight $w(\boldsymbol x)\equiv 1$.

\begin{theo}\label{theo:1-5}
A Euclidean $4$-design in $\mathbb R^n$ having the parameters
given in Theorem {\rm\ref{theo:1-4} (2) (ii)} exists if and only if a tight spherical
$4$-design on $S^n\subset\mathbb R^{n+1}$
exists.
\end{theo}

\noindent
{\bf Remark:} If $k=2$ and $k=3$ in the parameters
given above (Theorem \ref{theo:1-4} (2) (ii)),
then $n=5$ and $n=21$ respectively. 
The existence of spherical tight $4$-design on $S^5$ and $S^{21}$ are known. They are also known to be unique.
$S^{117}$, i.e. $k=6$, is the first case in which the existence of a spherical tight $4$-design is unknown (\cite{B-M-V}, see also \cite{B-B-5}).

\begin{theo}\label{theo:1-6}
\begin{enumerate}
\item The following is a family of feasible parameters
for tight Euclidean $4$-design in $\mathbb R^n$.

$n=(6k-3)^2-3$, with any positive integer $k$,

$|X_1|=(6k^2-6k+1)(36k^2-36k+7),\
|X_2|=3(36k^2-36k+7)(2k-1)^2$,

$A(X_1)=\left\{\frac{18k^2-27k+8}{6(9k^2-9k+1)(2k-1)},\ 
-\frac{18k^2-9k-1}{6(9k^2-9k+1)(2k-1)}
\right\}$,

$A(X_2)=\left\{\frac{36k^3-54k^2+25k-4}
{2(6k^2-6k+1)(18k^2-18k+5)},\
-\frac{36k^3-54k^2+25k-3}
{2(6k^2-6k+1)(18k^2-18k+5)}\right\}$,
 
 $A(X_1,X_2)=\left\{\sqrt{\frac{36k^2-36k+4}{(36k^2-36k+6)
 (36k^2-36k+10)}},\
 -\sqrt{\frac{36k^2-36k+10}{(36k^2-36k+6)
 (36k^2-36k+4)}}\right\}$,
 
 $r_1=1$,
 $r_2= \sqrt{\frac{3(18k^2-18k+5)(6k^2-6k+1)}
 {9k^2-9k+1}}
$, 

$w(\boldsymbol x)=1$ for $\boldsymbol x\in X_1$ and $w(\boldsymbol x)= \frac{1}{81(2k-1)^4}$
for $\boldsymbol x\in X_2$.

The intersection numbers of the
corresponding coherent configurations are given as polynomials of $k$ (see 
Appendix II).
\item If $2\leq n\leq 15^2-3$, then 
tight Euclidean $4$-design supported by 2 concentric
spheres is similar to one of the examples given in
Theorem I, Theorem II and Theorem III  in\cite{B-1} or to one of those having the parameters given above in this theorem. 
\end{enumerate}
\end{theo}

\noindent
{\bf Remark:} If $k=1$ in the parameters given above, then $n=6$ and the existence of the 
Euclidean tight $4$-design is known (Theorem I in\cite{B-1}). The first open parameters in this case is when $k=2$, i.e., $n=78$, which is also mentioned in \cite{B-1}.

In \S 2, we give some basic facts on the Euclidean $t$-designs.
In \S 3, we consider Euclidean $t$-designs having the structures
of coherent configurations and prove 
Theorem \ref{theo:1-2} and Theorem \ref{theo:1-3}.
In \S 4, we consider the Euclidean $4$-designs
supported by 2 concentric spheres and give the proof
for Theorem \ref{theo:1-4}, Theorem \ref{theo:1-5}
and Theorem \ref{theo:1-6}.

\section{Some basic facts on Euclidean $t$-designs}

As for the detailed definition and the basic properties
of Euclidean designs and examples of Euclidean designs please refer \cite{N-S,D-S,Baj-1,Baj-2,
B-B-3,B-B-4,B-B-5,B-B-S,B-B-H-S,B-1,B-2,M-0,M,
V-C}, etc. 
Here we only give the fact we need directly to prove our 
main theorems.
The following theorem gives a very useful
condition which is equivalent to
the definition of Euclidean $t$-designs. 
\begin{theo}[Neumaier-Seidel (see \cite{N-S})]
\label{theo:N-S}$ $\\
The following conditions are equivalent.
\begin{enumerate}
\item $(X,w)$ is a Euclidean $t$-design.
\item The following equation holds
$$\sum_{\boldsymbol x\in X}w(\boldsymbol x)
\|\boldsymbol x\|^{2j}\varphi_l(\boldsymbol x)=0 $$
for any harmonic polynomial
$\varphi_l\in\mbox{Harm}_l(\mathbb R^n)$, integers
$l$ and $j$ satisfying $1\leq l\leq t$ and 
$0\leq j\leq \frac{t-l}{2}$.
\item $$\sum_{\boldsymbol x\in X}w(\boldsymbol x)
f(\boldsymbol x)=\sum_{\boldsymbol x\in X}w(\boldsymbol x)
f(\tau(\boldsymbol x))$$
holds for any $f\in\mathcal P_t(\mathbb R^n)$ and
$\tau\in O(n)$, where $O(n)$ is the orthogonal group of degree $n$.
\end{enumerate}
\end{theo}
{\bf Remark:} Note that the condition (3) in 
Theorem \ref{theo:N-S} says that any kind of moments
of $X$ with degree at most $t$ is invariant
under any orthogonal transformations of 
$\mathbb R^n$. This concept is closely related
to the concept of rotatable designs
in statics (cf. \cite{B-H}).
Also note that Definition \ref{dfe:1-1}
is interpreted  as cubature formulas in analysis
(cf. Sobolev \cite{S-1,S-2} or \cite{M-0,M})
\\
 
Theorem \ref{theo:N-S} implies the following proposition
(see \cite{B-B-3}).
\begin{pro}
Let $(X,w)$ be a weighted finite set
in $\mathbb R^n$. Let $\rho$ be a similar transformation
of $\mathbb R^n$ fixing the origin. Let 
$\mu$ be a positive real number. 
Let $X'=\rho^{-1}(X)$, 
and $w'$ be a weight function on $X'$
defined by $w'(\boldsymbol  x')=\mu w(\rho(\boldsymbol x'))$
for any $\boldsymbol x'\in X'$.
The following conditions are equivalent.
\begin{enumerate}
\item $(X,w)$ is a Euclidean $t$-design.
\item $(X',w')$ is a Euclidean $t$-design.
\end{enumerate}
\label{pro:similar}
\end{pro}
We say that Euclidean $t$-designs $(X,w)$ and
 $(X',w')$ are similar if they satisfy the condition of 
 Proposition \ref{pro:similar}.
 Theorem \ref{theo:N-S} also implies
 the following.
 \begin{pro} Let $(X,w)$ be a weighted set in $\mathbb R^n$.
 Assume $0\not\in X$.
 Then $(X,w)$ is a Euclidean $t$-design if and only if
 $(X\cup\{0\},w)$ is a Euclidean $t$-design with $w(0)$ 
 any positive real number.
 \end{pro}
Let $h_l=h_{n,l}=\dim(\mbox{Harm}_l
({\mathbb R}^n))$ and
$\varphi_{l,1},\ \ldots\ ,\
\varphi_{l,h_l}$ be an
orthonormal basis of 
$\mbox{Harm}_l({\mathbb R}^n)$ 
with respect to the inner product
$\langle -,- \rangle$ defined by 
$$\langle \varphi,\psi \rangle=
\frac{1}{|S^{n-1}|}\int_{S^{n-1}}
\varphi(\boldsymbol x)\psi(\boldsymbol x)
d\sigma(\boldsymbol x) \quad 
\mbox{for}\ \varphi,\ \psi\in 
\mathcal P(\mathbb R^n).
$$
The following theorem is well known
(see \cite{E}).
\begin{theo}
Let $Q_l=Q_{n,l}$ be the Gegenbauer
polynomial of degree $l$ normalized 
so that satisfying $Q_l(1)=h_l$. Then
$$\sum_{i=1}^{h_l}\varphi_{l,i}(\boldsymbol x)
\varphi_{l,i}(\boldsymbol y)=Q_l(\boldsymbol x
\cdot \boldsymbol y)$$
holds for any $\boldsymbol x,\ \boldsymbol y\in
S^{n-1}$.
\end{theo}

Let $(X,w)$ be a Euclidean $t$-design
in $\mathbb R^n$.
Let $X=\cup_{i=1}^pX_i$.
and $r_i=\|\boldsymbol x\|$ for 
$\boldsymbol x\in X_i$, $1\leq i\leq p$.
For any nonnegative integers $l$ and $j$, we define matrices $H_{l,j}$  
whose rows 
and columns are indexed by $X$ and $\{\varphi_{l,1},\varphi_{l,2}
\ldots,\varphi_{l,h_l}\}$
respectively.
The $(\boldsymbol x, i)$-entry of $H_{l,j}$ for 
$\boldsymbol x\in
X_{\lambda}$
 is given by
$H_{l,j}(\boldsymbol  x,i)=\sqrt{w(\boldsymbol  x)}\|\boldsymbol  x\|^{2j}
\varphi_{l,i}(\boldsymbol  x)$.
Then the definition of Euclidean designs
implies the following proposition.
\begin{pro}\label{pro:3-2} Notation and definition are 
given as above. If $l_1+l_2+2j_1+2j_2\leq t$,
then the following holds. 
$$^tH_{l_1,j_1}\ H_{l_2,j_2}=
\left(\sum_{\lambda=1}^p
W(X_\lambda)
r_\lambda^{l_1+l_2+2(j_1+j_2)}\right)\Delta_{l_1,l_2},$$
where  $\Delta_{l_1,l_2}$ is the $0$ matrix
of size $h_{l_1}\times h_{l_2}$ for $l_1\not=l_2$ and
$\Delta_{l_1,l_1}$ is the identity matrix of size $h_{l_1}$.
\end{pro}

\section{Sufficient conditions for Euclidean designs
to have the structures
of coherent configurations}
Let $(X,w)$ be a Euclidean $t$-design 
supported by $p$ concentric spheres.
We use notation given in \S 1 and \S 2.
Let $\lambda, \mu$ be any integer satisfying
$1\leq \lambda, \mu\leq p$ and
$X_\lambda,X_\mu\not=\{0\}$.
For any $(\boldsymbol x,\boldsymbol y)\in X_\lambda\times
X_\mu$  let
$$p_{\alpha_{\lambda,\nu}^{(u)},\alpha_{\nu,\mu}^{(v)}}
(\boldsymbol x,\boldsymbol y)=
\left\{\boldsymbol z \in X_\nu\ \bigg|\
\frac{\boldsymbol x\cdot\boldsymbol z}{r_\lambda r_\nu}
=\alpha_{\lambda,\nu}^{(u)},\
\frac{\boldsymbol z\cdot\boldsymbol y}{r_\nu r_\mu}
=\alpha_{\nu,\mu}^{(v)} \right\}.$$
For any $(\boldsymbol x,\boldsymbol y)\in X_\lambda
\times X_\mu$ satisfying $\frac{\boldsymbol x\cdot\boldsymbol y}
{r_\lambda r_\mu}=\alpha_{\lambda,\mu}^{(q)}$,
the following holds.
\begin{equation}\label{equ:3-1}p_{\alpha_{\lambda,\lambda}^{(0)},\alpha_{\lambda,\mu}^{(u)}}
(\boldsymbol x,\boldsymbol y)=
p_{\alpha_{\lambda,\mu}^{(u)},\alpha_{\mu,\mu}^{(0)}}
(\boldsymbol x,\boldsymbol y)=\delta_{u,q}.
\end{equation}
We have the following proposition.

\begin{pro}\label{pro:3-1}
Let $(X,w)$ be a Euclidean
$t$-design supported by $p$ concentric spheres. Assume that $0\not \in X$ and the weight function 
is constant on each $X_\nu (1\leq \nu\leq p)$,
i.e., $w(\boldsymbol x)\equiv w_\nu$
for any
$\boldsymbol x\in X_\nu (1\leq \nu\leq p)$.
Then the followings hold for any nonnegative integers $l,\ k$ and $j$ satisfying $l+k+2j\leq t$.
\begin{enumerate}
\item For $\boldsymbol x, \boldsymbol y\in
X_\lambda$ and $\frac{\boldsymbol x\cdot \boldsymbol y}
{r_\lambda^2}
=\alpha_{\lambda,\lambda}^{(q)}$,
\begin{eqnarray}
&&
\sum_{\nu=1}^p\sum_{u=1}^{s_{\lambda,\nu}}
\sum_{v=1}^{s_{\nu,\lambda}}
w_\nu r_\nu^{l+k+2j}
Q_l(\alpha_{\lambda,\nu}^{(u)})Q_k(\alpha_{\nu,\lambda}^{(v)})
p_{\alpha_{\lambda,\nu}^{(u)},\alpha_{\nu,\lambda}^{(v)}}(
\boldsymbol x,\boldsymbol y)
\nonumber\\
&&=\delta_{l,k}Q_l(\alpha_{\lambda,\lambda}^{(q)})
\sum_{\nu=1}^p|X_\nu| w_\nu
r_\nu^{2l+2j}
\nonumber\\
&&-
\left\{\begin{array}{ll}w_\lambda r_\lambda^{l+k+2j}
\left(Q_l(\alpha_{\lambda,\lambda}^{(q)})Q_k(1)
+Q_l(1)Q_k(\alpha_{\lambda,\lambda}^{(q)})\right)
&\mbox{for}\ q\neq 0\\
w_\lambda r_\lambda^{l+k+2j}Q_l(1)Q_k(1)
&\mbox{for}\ q=0.
\end{array}
\right.
\label{equ:3-2}
\end{eqnarray}
\item For $\boldsymbol x\in X_\lambda$,
$\boldsymbol y\in X_\mu$, $\lambda\neq\mu$
and $\frac{\boldsymbol x\cdot\boldsymbol y}
{{r_\lambda}{ r_\mu}}
=\alpha_{\lambda,\mu}^{(q)}$,
\begin{eqnarray}
&&
\sum_{\nu=1}^p\sum_{u=1}^{s_{\lambda,\nu}}
\sum_{v=1}^{s_{\nu,\mu}}
w_\nu r_\nu^{l+k+2j}
Q_l(\alpha_{\lambda,\nu}^{(u)})Q_k(\alpha_{\nu,\mu}^{(v)})
p_{\alpha_{\lambda,\nu}^{(u)},\alpha_{\nu,\mu}^{(v)}}
(\boldsymbol x,\boldsymbol y)
\nonumber\\
&&=\delta_{l,k}Q_l(\alpha_{\lambda,\mu}^{(q)})
\sum_{\nu=1}^p|X_\nu| w_\nu
r_\nu^{2l+2j}-
w_\lambda r_\lambda^{l+k+2j}
Q_l(1)Q_k(\alpha_{\lambda,\mu}^{(q)})
\nonumber\\
&&-w_\mu r_\mu^{l+k+2j}
Q_l(\alpha_{\lambda,\mu}^{(q)})Q_k(1).
\label{equ:3-3}
\end{eqnarray}
\end{enumerate}
\end{pro}
{\bf Proof} Choose non negative integers
$j_1$ and $j_2$ satisfying $j_1+j_2=j$. Then 
Proposition \ref{pro:3-2}
implies
\begin{equation}(H_{l,j_1}\ {^tH_{l,j_1}})(H_{k,j_2}\ {^tH_{k,j_2}})
=\delta_{l,k}
\sum_{\nu=1}^p|X_\nu| w_\nu
r_\nu^{2j_1+2j_2+2l}(H_{l,j_1}\ {^tH_{l,j_2}}).
\label{equ:3-4}
\end{equation}
The
$(\boldsymbol x,\boldsymbol y)$-entry
of the left hand side of (\ref{equ:3-4})
gives
\begin{eqnarray}
&&((H_{l,j_1}\ {^tH_{l,j_1}})(H_{k,j_2}\ {^tH_{k,j_2}}))
(\boldsymbol x,\boldsymbol y)
=
\sum_{\boldsymbol z\in X}(H_{l,j_1}\ {^tH_{l,j_1}})
(\boldsymbol x,\boldsymbol z)
(H_{k,j_2}\ {^tH_{k,j_2}})(\boldsymbol z,\boldsymbol y)
\nonumber\\
&&=\sqrt{w(\boldsymbol x)w(\boldsymbol y)}\|\boldsymbol x\|^{l+2j_1}\|\boldsymbol y\|^{k+2j_2}\sum_{\boldsymbol z\in X}
w(\boldsymbol z)\|\boldsymbol z\|^{l+k+2j}
Q_l\left(\frac{\boldsymbol x\cdot
\boldsymbol z}
{\|\boldsymbol x\|\|\boldsymbol z\|}\right)
Q_k\left(\frac{\boldsymbol z\cdot\boldsymbol y}
{\|\boldsymbol z\|\|\boldsymbol y\|}\right).
\nonumber\\
&&=\sqrt{w_\lambda w_\mu}
 r_\lambda^{l+2j_1}r_\mu^{k+2j_2}\nonumber\\
&&
\times\left(\sum_{\nu=1}^p
w_\nu r_\nu^{l+k+2j}
\sum_{u=1-\delta_{\nu,\lambda}}^{s_{\lambda,\nu}}\
\sum_{v=1-\delta_{\nu,\mu}}^{s_{\nu,\mu}}
p_{\alpha_{\lambda,\nu}^{(u)},\alpha_{\nu,\mu}^{(v)}}
(\boldsymbol x,\boldsymbol y)
Q_l(\alpha_{\lambda,\nu}^{(u)})
Q_k(\alpha_{\nu,\mu}^{(v)})\right).
\label{equ:3-5}
\end{eqnarray}
On the other hand, the right hand side of (\ref{equ:3-4}) gives
\begin{eqnarray}
&&\delta_{l,k}
\sum_{\nu=1}^p|X_\nu| w_\nu
r_\nu^{2j_1+2j_2+2l}
\sqrt{w_\lambda w_\mu}
r_\lambda^{l+2j_1}r_\mu^{l+2j_2}
Q_l\left(\frac{\boldsymbol x\cdot \boldsymbol y}
{\|\boldsymbol x \|\|\boldsymbol y\|}\right).
\label{equ:3-6}
\end{eqnarray}
Since $\sqrt{w_\lambda w_\mu}r_\lambda^{l+2j_1}
r_\mu^{l+2j_2}\neq0$ and $j=j_1+j_2$, 
(\ref{equ:3-1}),
(\ref{equ:3-4}), (\ref{equ:3-5}) and (\ref{equ:3-6}) imply
(\ref{equ:3-2}) and (\ref{equ:3-3}).
\hfill\qed\\

\noindent
{\bf Proof of Theorem \ref{theo:1-2} with the condition (1)}\\
If $0\in X$, then $\varepsilon_S=1$
and $X\backslash \{0\}$ is a Euclidean
$t$ design on the union of $p-1$ concentric spheres $S'=
S\backslash \{0\}$
and $s_{\lambda,\nu}+s_{\nu,\mu}\leq
t-2((p-1)-\varepsilon_{S'}-2)$ holds.
It is easy to see that 
if $X\backslash \{0\}$ has the structure of a coherent configuration, then $ X$ also has the structure of a coherent configuration.
Therefore in the following we assume $0\not \in X$ 
($\varepsilon_S=0$).
For each fixed $l,k$, both (\ref{equ:3-2}) and (\ref{equ:3-3}) 
consist of $[\frac{t-l-k}{2}]+1$ linear equations of 
indeterminate 
$p_{\alpha_{\lambda,\nu}^{(u)},\alpha_{\nu,\mu}^{(v)}}
(\boldsymbol x,\boldsymbol y)$ ($j=0,1,\ldots,[\frac{t-l-k}{2}]$).
For each $l,k,j$, right hand sides of the both linear equations 
(\ref{equ:3-2}) and (\ref{equ:3-3}) are functions of 
$\alpha_{\lambda,\mu}^{(q)}$,
say $F_{l,k,j}(\alpha_{\lambda,\mu}^{(q)})$ and
independent of the choice of 
$ \boldsymbol x\in X_\lambda$
and $\boldsymbol y\in X_\mu$
whenever 
$\boldsymbol x\cdot\boldsymbol y={r_\lambda}
{r_\mu}
\alpha_{\lambda,\mu}^{(q)}$ is satisfied.
Let us consider the left hand sides of 
(\ref{equ:3-2}) and (\ref{equ:3-3}) together.
To do so we consider the following 
system of linear equations.
\begin{equation}
\sum_{\nu=1}^p
\sum_{u=1}^{s_{\lambda,\nu}}
\sum_{v=1}^{s_{\nu,\mu}}
w_\nu r_\nu^{l+k+2j}
Q_l(\alpha_{\lambda,\nu}^{(u)})Q_k(\alpha_{\nu,\mu}^{(v)})
p_{\alpha_{\lambda,\nu}^{(u)},\alpha_{\nu,\mu}^{(v)}}
(\boldsymbol x,\boldsymbol y)
=F_{l,k,j}(\alpha_{\lambda,\mu}^{(q)})
\end{equation}
Let 
$\Psi_{l,k,\nu}(\boldsymbol x,\boldsymbol y)=w_\nu r_\nu^{l+k}\sum_{u=1}^{s_{\lambda,\nu}}
\sum_{v=1}^{s_{\nu,\mu}}
Q_l(\alpha_{\lambda,\nu}^{(u)})Q_k(\alpha_{\nu,\mu}^{(v)})
p_{\alpha_{\lambda,\nu}^{(u)},\alpha_{\nu,\mu}^{(v)}}
(\boldsymbol x,\boldsymbol y)$.
Then for any non negative pair $(l,k)$ of integers,
satisfying $p-1\leq \frac{t-l-k}{2}$, 
we obtain system of the following $p$ equations
with indeterminates $\{\Psi_{l,k,\nu}(\boldsymbol x,\boldsymbol y)
\mid
\nu=1,\ldots,p\}$.
\begin{equation}\sum_{\nu=1}^pr_\nu^{2j}
\Psi_{l,k,\nu}(\boldsymbol x,\boldsymbol y)=
F_{l,k,j}(\alpha_{\lambda,\mu}^{(q)}),
\quad j=0,1,\ldots,p-1.
\label{equ:1-8}
\end{equation}
Since the coefficient matrix of the linear 
equations (\ref{equ:1-8}) equals
$$\left[
\begin{array}{ccccc}
1&\cdots&1&\cdots&1\\
r_1^2&\cdots&r_\nu^2&\cdots&r_p^2\\
\vdots&\cdots&\cdots&\cdots&\vdots\\
r_1^{2(p-1)}&\cdots&r_\nu^{2(p-1)}&\cdots&r_p^{2(p-1)}
\end{array}
\right]
$$
which is invertible.
Hence, for each non negative pair $(l,k)$ of integers, satisfying
$l+k\leq t-2p+2$, and $\nu$, $1\leq \nu\leq p$, $\Psi_{l,k,\nu}(\boldsymbol x,\boldsymbol y)$ determined
uniquely by $\alpha_{\lambda,\mu}^{(q)}$
independent of the choice of $(\boldsymbol x,\boldsymbol y)
\in X_\lambda\times X_\mu$ satisfying
$\frac{\boldsymbol x\cdot\boldsymbol y}{r_\lambda r_\mu}=
\alpha_{\lambda,\mu}^{(q)}$.
More precisely, for each non negative pair $(l,k)$ of integers
satisfying $l+k\leq  t-2(p-1)$, and $\nu$ with 
$1\leq \nu\leq p$, 
\begin{equation}\sum_{u=1}^{s_{\lambda,\nu}}
\sum_{v=1}^{s_{\nu,\mu}} 
Q_l(\alpha_{\lambda,\nu}^{(u)})Q_k(\alpha_{\nu,\mu}^{(v)})
p_{\alpha_{\lambda,\nu}^{(u)},\alpha_{\nu,\mu}^{(v)}}
(\boldsymbol x,\boldsymbol y)=
\frac{G_{l,k,\nu}(\alpha_{\lambda,\mu}^{(q)})}{w_\nu r_\nu^{l+k}}
\label{equ:1-9}
\end{equation}
holds, 
 where
$G_{l,k,\nu}( \alpha_{\lambda,\mu}^{(q)})$
depends only on $l$, $k$, $\nu$ and $\alpha_{\lambda,\mu}^{(q)}$.
Since $s_{\lambda,\nu}+s_{\nu,\mu}\leq t-2(p-2)$,
then
(\ref{equ:1-9}) holds for any
$l$ and $k$ satisfying
$0\leq l\leq s_{\lambda,\nu}-1$ and 
$0\leq k\leq s_{\nu,\mu}-1$. 
Then (\ref{equ:1-9}) gives a
system of linear equations whose
coefficient matrices are
the tensor product
 $$\left[
\begin{array}{ccccc}
1&\cdots&1&\cdots&1\\
Q_1(\alpha_{\lambda,\nu}^{(1)})&\cdots
&Q_1(\alpha_{\lambda,\nu}^{(u)})
&\cdots&Q_1(\alpha_{\lambda,\nu}^{(s_{\lambda,\nu})})\\
\vdots&\cdots&\cdots&\cdots&\vdots\\
Q_{s_{\lambda,\nu}-1}(\alpha_{\lambda,\nu}^{(1)})
&\cdots&Q_{s_{\lambda,\nu}-1}(\alpha_{\lambda,\nu}^{(u)})
&\cdots&Q_{s_{\lambda,\nu}-1}
(\alpha_{\lambda,\nu}^{(s_{\lambda,\nu})})
\end{array}
\right]$$
$$\bigotimes
\left[
\begin{array}{ccccc}
1&\cdots&1&\cdots&1\\
Q_1(\alpha_{\nu,\mu}^{(1)})&\cdots
&Q_1(\alpha_{\nu,\mu}^{(u)})
&\cdots&Q_1(\alpha_{\nu,\mu}^{(s_{\nu,\mu})})\\
\vdots&\cdots&\cdots&\cdots&\vdots\\
Q_{s_{\nu,\mu}-1}(\alpha_{\nu,\mu}^{(1)})
&\cdots&Q_{s_{\nu,\mu}-1}(\alpha_{\nu,\mu}^{(u)})
&\cdots&Q_{s_{\nu,\mu}-1}
(\alpha_{\nu,\mu}^{(s_{\nu,\mu})})
\end{array}
\right]
$$
of two invertible matrices.
Hence $p_{\alpha_{\lambda,\nu}^{(u)},\alpha_{\nu,\mu}^{(v)}}
(\boldsymbol x,\boldsymbol y)$
determined uniquely by
$\alpha_{\lambda,\nu}^{(u)},\ \alpha_{\nu,\mu}^{(v)}$
and $\alpha_{\lambda,\mu}^{(q)}$ which does not depend of the
choice of $(\boldsymbol x,\boldsymbol y)\in 
X_\lambda\times X_\mu$
satisfying $\frac{\boldsymbol x\cdot\boldsymbol y}
{{r_\lambda}{r_\mu}}
=\alpha_{\lambda,\mu}^{(q)}$.
This completes the proof of Theorem \ref{theo:1-2} with the condition (1).\\

Next, we consider the case when $X$ is antipodal. 
Let $\lambda$ and $\mu$ be any integers
satisfying $1\leq \lambda,\mu\leq p$ and $X_\lambda,X_\mu\not=\{0\}$.
Since $X$ is antipodal, $-1\in A(X_\lambda,X_\lambda)$ holds.
Let us denote $\alpha_{\lambda,\lambda}^{(1)}=-1$.
Also $-\alpha_{\lambda,\mu}^{(u)}\in A(X_\lambda,X_\mu)$
for any $u$ satisfying $1\leq u\leq s_{\lambda,\mu}$.
If $-\alpha_{\lambda,\mu}^{(u)}\in A(X_\lambda,X_\mu)$, then 
let $\alpha_{\lambda,\mu}^{(u^*)}
=-\alpha_{\lambda,\mu}^{(u)}$.
For any $(\boldsymbol x,\boldsymbol y)\in X_\lambda
\times X_\mu$ satisfying 
$\frac{\boldsymbol x\cdot\boldsymbol y}{r_\lambda r_\mu}
=\alpha_{\lambda,\mu}^{(q)}$, the following holds. 

\begin{equation}p_{\alpha_{\lambda,\lambda}^{(1)}
,\alpha_{\lambda,\mu}^{(u)}}
(\boldsymbol x,\boldsymbol y)=p_{\alpha_{\lambda,\mu}^{(u)},\alpha_{\mu,\mu}^{(1)}}(\boldsymbol x,\boldsymbol y)
=\delta_{u, q^*}.
\label{equ:1-11}
\end{equation} 
Then similar arguments as before give the following
proposition.
\begin{pro} Let $X$ be a Euclidean 
$t$-design. Assume $X$ is antipodal, $0\not\in X$ 
and the weight function 
is constant on each $X_\nu (1\leq \nu\leq p)$,
i.e., $w(\boldsymbol x)\equiv w_\nu$
for any
$\boldsymbol x\in X_\nu (1\leq \nu\leq p)$.
Then the following hold 
for any nonnegative integers 
$l, k$ and $j$ satisfying $l+k+2j\leq t$.
\begin{enumerate}
\item For $\boldsymbol x, \boldsymbol y\in
X_\lambda$ and $\frac{\boldsymbol x\cdot \boldsymbol y}
{r_\lambda^2}
=\alpha_{\lambda,\lambda}^{(q)}$,
\begin{eqnarray}
&&\sum_{\nu=1}^p
w_\nu r_\nu^{l+k+2j}
\sum_{u=1+\delta_{\lambda,\nu}}^{s_{\lambda,\nu}}\
\sum_{v=1+\delta_{\nu,\lambda}}^{s_{\nu,\lambda}}
p_{\alpha_{\lambda,\nu}^{(u)},\alpha_{\nu,\lambda}^{(v)}}
(\boldsymbol x,\boldsymbol y)
Q_l(\alpha_{\lambda,\nu}^{(u)})
Q_k(\alpha_{\nu,\lambda}^{(v)})
\nonumber\\
&&
\nonumber\\
&&=\delta_{l,k}Q_l(\alpha_{\lambda,\lambda}^{(q)})
\sum_{\nu=1}^p|X_\nu| w_\nu
r_\nu^{2l+2j}
\nonumber\\
&&-
\left\{\begin{array}{ll}
w_\lambda r_\lambda^{l+k+2j}
((-1)^{l+k}+1)\left(Q_l(\alpha_{\lambda,\lambda}^{(q)})Q_k(1)
+Q_l(1)Q_k(\alpha_{\lambda,\lambda}^{(q)})\right)
&\mbox{for}\ q\neq 0,1\\
((-1)^{l+k}+1)w_\lambda r_\lambda^{l+k+2j}
Q_l(1)Q_k(1)&\mbox{for}\ q=0\\
\left((-1)^k+(-1)^l\right)
w_\lambda r_\lambda^{l+k+2j}Q_l(1)Q_k(1)
&\mbox{for}\ q=1.
\end{array}
\right.
\nonumber\\
&&
\label{equ:1-12}
\end{eqnarray}
\item For $\boldsymbol x\in X_\lambda$,
$\boldsymbol y\in X_\mu$. Assume $\lambda\neq\mu$
and $\frac{\boldsymbol x\cdot\boldsymbol y}
{r_\lambda r_\mu}
=\alpha_{\lambda,\mu}
^{(q)}$.
Then 
\begin{eqnarray}
&&
\sum_{\nu=1}^p
w_\nu r_\nu^{l+k+2j}
\sum_{u=1+\delta_{\lambda,\nu}}^{s_{\lambda,\nu}}\
\sum_{v=1+\delta_{\nu,\mu}}^{s_{\nu,\mu}}
p_{\alpha_{\lambda,\nu}^{(u)},\alpha_{\nu,\mu}^{(v)}}(
\boldsymbol x,\boldsymbol y)
Q_l(\alpha_{\lambda,\nu}^{(u)})
Q_k(\alpha_{\nu,\mu}^{(v)})
\nonumber\\
&&
\nonumber\\
&&=\delta_{l,k}Q_l(\alpha_{\lambda,\mu}^{(q)})
\sum_{\nu=1}^p|X_\nu| w_\nu
r_\nu^{2l+2j}\nonumber\\
&&
-\left((-1)^{l+k}+1\right)
\left(w_\mu r_\mu^{l+k+2j}
Q_l(\alpha_{\lambda,\mu}^{(q)})Q_k(1)
+w_\lambda r_\lambda^{l+k+2j}
Q_l(1)Q_k(\alpha_{\lambda,\mu}^{(q)})\right)\nonumber\\
\label{equ:1-13}
\end{eqnarray}
\end{enumerate}
\end{pro}
{\bf Proof} (\ref{equ:1-11}) and (\ref{equ:3-4}),
(\ref{equ:3-5}), (\ref{equ:3-6}) imply 
(\ref{equ:1-12}) and (\ref{equ:1-13}).
\hfill\qed\\

\noindent
{\bf Proof of Theorem \ref{theo:1-2} with the condition (2)}\\
 As we mentioned at the beginning of the proof for
Theorem \ref{theo:1-2} with the condition (1), 
it is enough if we prove
the statement for the case $0\not\in X$, i.e. $\varepsilon_S=0$.  The same argument as we used
in the proof of Theorem \ref{theo:1-2} with the condition (1) implies
that for any 
$(\boldsymbol x,\boldsymbol y)\in X_\lambda\times X_\mu$ satisfying $\frac{\boldsymbol x\cdot\boldsymbol y}
{r_\lambda r_\mu}=\alpha_{\lambda,\mu}^{(q)}$
the following holds
\begin{equation}
\sum_{u=1+\delta_{\lambda,\nu}}^{s_{\lambda,\nu}}\
\sum_{v=1+\delta_{\nu,\mu}}^{s_{\nu,\mu}}
p_{\alpha_{\lambda,\nu}^{(u)},\alpha_{\nu,\mu}^{(v)}}
(\boldsymbol x,\boldsymbol y)
Q_l(\alpha_{\lambda,\nu}^{(u)})
Q_k(\alpha_{\nu,\mu}^{(v)})
=\frac{G_{l,k,\nu}(\alpha_{\lambda,\nu}^{(q)})}{w_\nu r_\nu^{l+k}}
\label{equ:1-14}
\end{equation}
for any positive integer $\nu$
and non negative integers 
$l,k$ satisfying $1\leq \nu\leq p$ and
$0\leq l+k\leq t-2(p-1)$,
where $G_{l,k,\nu}(\alpha_{\lambda,\nu}^{(q)})$
is independent of the choice of $(\boldsymbol x,\boldsymbol y)$.
Since $s_{\lambda,\nu}+s_{\nu,\mu}-\delta_{\lambda,\nu}
-\delta_{\nu,\mu}\leq t-2(p-2)$,
(\ref{equ:1-14}) holds for any
$l$ and $k$ satisfying
$0\leq l\leq s_{\lambda,\nu}-\delta_{\lambda,\nu}-1$
and 
$0\leq k\leq s_{\nu,\mu}-\delta_{\nu,\mu}-1$.
Then for each triple
$\lambda,\nu,\mu$ we obtain a system of linear equations
with determinates $p_{\alpha_{\lambda,\nu}^{(u)},\alpha_{\nu,\mu}^{(v)}}
(\boldsymbol x,\boldsymbol y)$
whose coefficient matrix is nonsingular.
This implies that the intersection numbers
$p_{\alpha_{\lambda,\nu}^{(u)},\alpha_{\nu,\mu}^{(v)}}
(\boldsymbol x,\boldsymbol y)$ is independent
of the choice of $\boldsymbol x\in X_\lambda,
\boldsymbol y\in X_\mu$
satisfying  $\frac{\boldsymbol x\cdot
\boldsymbol y}{r_\lambda r_\mu}=
\alpha_{\lambda,\mu}^{(q)}$.
This completes the proof of Theorem \ref{theo:1-2} with the condition (2). 
\hfill\qed\\

\noindent
{\bf Proof of Theorem \ref{theo:1-3}}\\
If $t=2e$ or $t=2e+1$, and $0\in X$, 
then $e$ must be an even integer
and $\frac{e}{2}+1=2(=p)$ (Proposition 2.4.4 and 2.4.5
in \cite{B-B-H-S}). Hence $t=4$ 
or $t=5$ and 
$X\backslash\{0\}$ is a tight spherical design 
having the structure of Q-polynomial association 
scheme.
Assume $0\not\in X$. 
If $t=2e$, then the arguments in the proof for Lemma 1.10 in  
\cite{B-B-3} imply that $w$ is constant on each $X_\lambda$ and $s_{\lambda,\mu}\leq e$ for 
any $1\leq\lambda,\mu\leq 2(=p)$. On the other hand,
if $t=2e+1$, then Proposition 2.4.6 in \cite{B-B-H-S}
and the arguments in the proof for
Lemma 1.7 in \cite{B-1} imply
that $X$ is antipodal, weight function is constant on
each $X_\lambda$, $s_{\lambda,\lambda}\leq e+1$,
$s_{\lambda,\mu}\leq e$ for any $1\leq \lambda\neq \mu\leq 2(=p)$.
Hence $s_{\lambda,\nu}-
\delta_{\lambda,\nu}=e$ holds for any $1\leq \lambda,\nu\leq 2
(=p)$. 
This implies
$s_{\lambda,\nu}-
\delta_{\lambda,\nu}+s_{\nu,\mu}-
\delta_{\nu,\mu}\leq 2e<t-2(p-2)=2e+1$.
If $t=1$, then $X=\{\boldsymbol x,\ -\boldsymbol x\}$ 
and it is on a sphere in $\mathbb R^n$ and $p=1$.
Hence Theorem \ref{theo:1-2} implies Theorem \ref{theo:1-3}.
\hfill\qed\\

\section{Euclidean $4$-designs on $2$ concentric
spheres and coherent configurations}  

In this section we consider a Euclidean $4$-design
$(X,w)$ supported by 2 concentric spheres.
We assume that $0 \not \in X$ and the weight function $w$ is constant on each
layer $X_1$ and $X_2$.
If $s_{\lambda,\mu}\leq 2$ for any 
$\lambda,\mu\in\{1,2\}$, then Theorem \ref{theo:1-2}
implies that $X$ has the structure of a
coherent configuration.\\
 
 \subsection{Proof for Theorem \ref{theo:1-4}}
\noindent
{\bf Proof of Theorem \ref{theo:1-4} (1)}\\
Theorem 2.3 in \cite{B-1} implies that
both $X_1$ and $X_2$ are spherical $2$-designs.
 If $s_{1,2}=1$ and $A(X_1,X_2)=\{\gamma\}$,
then $X_2\subset \{\boldsymbol x
\mid \boldsymbol x\cdot\boldsymbol u=r_2\gamma\}$
where $\boldsymbol u$ is any fixed point in $X_1$.
Thus $X_2$ is on the intersection of the two
$(n-1)$-dimensional spheres, 
$S_2$ and the sphere $\{\boldsymbol x
\mid \boldsymbol x\cdot\boldsymbol u=r_2\gamma\}$
centered at $\boldsymbol u$. Hence $X_2$ is on an
$n-2$ dimensional sphere and 
$X_2$ cannot be a spherical 2-design on   
$(n-1)$-dimensional sphere $S_2$.
Hence we must have $s_{1,2}=2$.\hfill\qed\\

Let $N_i=|X_i|$ for $i=1,2$. 
By Proposition \ref{pro:similar} (in \S 2 of this paper)
we may assume the following:

$N_2\geq N_1$ and $r_1=1$, $w(\boldsymbol x)
\equiv 1$ on $X_1$ and $w(\boldsymbol x)
\equiv w_2$ on $X_2$.
Since $(X,w)$ is a Euclidean 4-design and  $X_i$, $i=1,2$ is a spherical $2$-design
we must have $|X|=N_1+N_2\geq {n+2\choose 2}$ 
and $N_1\geq n+1$.\\

We first prove the following theorem. 
\begin{theo}
\label{theo:4-1} 
Let $(X,w)$ be a Euclidean $4$-design 
on $2$ concentric spheres in $\mathbb R^n$. 
If $N_1= n+1$, then $(X,w)$ is a Euclidean
tight $4$-design.
\end{theo}

\noindent
{\bf Proof}\quad 
Since $N_1= n+1$, $X_1$ is a tight spherical $2$-design.
Hence $X_1$ is a regular simplex, i.e. $s_{1,1}=1$, 
on the unit 
sphere $S_1=S^{n-1}$ (see \cite{D-G-S,B-B-1}).  
If $N_1+N_2= {n+2\choose 2}$, then $(X,w)$ is a Euclidean
tight $4$-design.
Hence we may assume $N_2\geq \frac{n(n+1)}{2}+1\geq n+2$.
Hence $X_2$ must be a $2$-distance set on a sphere, that is,
$s_{2,2}=2$ and $N_2\leq \frac{n(n+3)}{2}$ holds (see \cite{D-G-S,B-B-1}).  
 
Let $p_{\alpha_{\lambda,\nu}^{(u)},\alpha_{\nu,\mu}^{(v)}}
(\boldsymbol x,\boldsymbol y)=p_{\alpha_{\lambda,\nu}^{(u)},\alpha_{\nu,\mu}^{(v)}}^{\alpha_{\lambda,\mu}^{(q)}}$
for $(\boldsymbol x,\boldsymbol y)\in X_\lambda\times X_\mu$
satisfying 
$\frac{\boldsymbol x\cdot\boldsymbol y}{r_\lambda,r_\mu}
=\alpha_{\lambda,\mu}^{(q)}$.
Let $\alpha_i=\alpha_{1,1}^{(i)}$ for $i=0,1$, 
$\beta_i=\alpha_{2,2}^{(i)}$ for $i=0,1,2$, and
$\gamma_i=\alpha_{1,2}^{(i)}$ for $i=1,2$.
We note that $\alpha_1=-\frac{1}{n}$ holds.
We assume $\beta_1>\beta_2$ and $\gamma_1>\gamma_2$.
Using the equations given in Proposition \ref{pro:3-1},
we can determine intersection numbers
of the corresponding coherent configuration.
By definition, we have the following immediately.

$p_{\gamma_1,\gamma_2}^{\alpha_0}=p_{\gamma_2,\gamma_1}^{\alpha_0}=p_{\gamma_1,\gamma_2}^{\beta_0}
=p_{\gamma_2,\gamma_1}^{\beta_0}
=p_{\beta_1,\beta_2}^{\beta_0}
=p_{\beta_2,\beta_1}^{\beta_0}=0$,

$p_{\alpha_1,\alpha_1}^{\alpha_0}=n$,\quad
$p_{\alpha_1,\alpha_1}^{\alpha_1}=n-1$,

$p_{\gamma_2,\gamma_2}^{\alpha_0}
=N_2-p_{\gamma_1,\gamma_1}^{\alpha_0}$,
\quad
$p_{\gamma_2,\gamma_2}^{\beta_0}=
n+1-p_{\gamma_1,\gamma_1}^{\beta_0}$,

$p_{\beta_2,\beta_2}^{\beta_0}
= N_2-1-p_{\beta_1,\beta_1}^{\beta_0}
$.\\
We also have the following.

$p_{\gamma_1,\gamma_2}^{\alpha_1}
=p_{\gamma_2,\gamma_1}^{\alpha_1}$,\quad

$p_{\beta_1,\beta_2}^{\beta_i}=p_{\beta_2,\beta_1}^{\beta_i}$, 
$p_{\gamma_1,\gamma_2}^{\beta_i}=p_{\gamma_2,\gamma_1}^{\beta_i}$
for $i=1,2$,

$p_{\alpha_i,\beta_j}^{\gamma_k}=
p_{\beta_j,\alpha_i}^{\gamma_k}
=0$ for $i=0,1$, $j=0,1,2$,
and $k=1,2$.

$p_{\gamma_i,\alpha_0}^{\gamma_k}
=p_{\alpha_0,\gamma_i}^{\gamma_k}
=p_{\gamma_i,\beta_0}^{\gamma_k}
=p_{\beta_0,\gamma_i}^{\gamma_k}\delta_{i,k}$,
$p_{\gamma_i,\alpha_1}^{\gamma_k}
=p_{\alpha_1,\gamma_i}^{\gamma_k}$,
$p_{\gamma_i,\beta_j}^{\gamma_k}
=p_{\beta_j,\gamma_i}^{\gamma_k}$ 
for $i,j,k=1,2$.\\
Proposition \ref{pro:3-1} (1)
with $\lambda=1$, $q=0$, $k=j=0$
for $l=1$ and $l=2$ imply
\begin{eqnarray}&&p_{\gamma_1,\gamma_1}^{\alpha_0}
=\frac{ N_2}{n\gamma_1^2+1},\label{equ:n+1-m2}\\
&&
\gamma_1\gamma_2 = -\frac{1}{n}.
\label{equ:gamma}
\end{eqnarray}
Therefore $\gamma_1>0>\gamma_2$ holds.
Next Proposition \ref{pro:3-1} (1)
with $\lambda=2$, $q=0$, $l=1$, $k=j=0$
and
$\lambda=2$, $q=0$, $l=1$, $k=0$, $j=1$ imply
\begin{eqnarray}&&p_{\beta_1,\beta_1}^{\beta_0}=
-\frac{\beta_2(N_2-1)+1}{\beta_1-\beta_2},\\
&&p_{\gamma_1,\gamma_1}^{\beta_0}
=\frac{n+1}{n\gamma_1^2+1}\label{equ:n+1-m1}.
\end{eqnarray}
Proposition \ref{pro:3-1} (1)
with $\lambda=2$, $q=0$, $l=2$, $k=0$,
$j=1$ inplies
$$n(\beta_1(N_2-1)+1)\beta_2+N_2-
n(1-\beta_1)=0.
$$
If $(\beta_1(N_2-1)+1)=0$ holds, then we must have
$N_2-n(1-\beta_1)=0$. This implies
$N_2=n+1$. This contradict our assumption $N_2>n+1$.
Hence $(\beta_1(N_2-1)+1)\not=0$
and we obtain 
\begin{eqnarray}&&\beta_2= -\frac{N_2+n\beta_1-n}{n(1+N_2\beta_1-\beta_1)}
\label{equ:n+1-beta2}\\
&&p_{\beta_1,\beta_1}^{\beta_0}
= \frac{N_2(N_2-n-1)}
{n(N_2-1)\beta_1^2+2n\beta_1+N_2-n}
\end{eqnarray}
Proposition \ref{pro:3-1} (1) with
$\lambda=1$, $q=1$, $(l,k,j)=(0,0,0),
(1,0,0)$ and $(1,1,0)$ imply
\begin{equation}p_{\gamma_1,\gamma_1}^{\alpha_1}
= \frac{N_2(1-\gamma_1^2)}{(\gamma_1^2n+1)^2},
\quad
p_{\gamma_2,\gamma_2}^{\alpha_1}
=\frac{ N_2\gamma_1^2(n^2\gamma_1^2-1)}{(\gamma_1^2n+1)^2},
\quad
p_{\gamma_1,\gamma_2}^{\alpha_1} =
 \frac{(n+1)N_2\gamma_1^2}{(\gamma_1^2n+1)^2}.
 \end{equation}
 Since $p_{\gamma_2,\gamma_2}^{\alpha_1}\geq 0$,
 we must have $0<\gamma_1\leq \frac{1}{\sqrt{n}}$.\\
Proposition \ref{pro:3-1} (1) with
$\lambda=2$, $q=1$, $(l,k,j)=(0,0,0),
(0,0,1), (1,0,0), (1,0,1)$,
$(1,1,0),$ and 
$(1,1,1)$ imply

$p_{\beta_1,\beta_1}^{\beta_1}
=\frac{N_2\bigg(n(N_2-1)(N_2-2n-1)\beta_1^3
-3\beta_1^2n^2-3n\beta_1+(N_2-n-2)(N_2-n)
\bigg)}{(2n\beta_1+n\beta_1^2N_2-n\beta_1^2+N_2-n)^2},
$

$p_{\beta_2,\beta_2}^{\beta_1}
= \frac{nN_2\beta_1(1+N_2\beta_1-\beta_1)^2(n\beta_1+1)}{(n(N_2-1)\beta_1^2+2n\beta_1
+N_2-n)^2},
$

$p_{\beta_1,\beta_2}^{\beta_1}=
p_{\beta_2,\beta_1}^{\beta_1}=
 \frac{n(-\beta_1+1)(N_2+n\beta_1-n)(1+N_2\beta_1-\beta_1)^2}{(2n\beta_1+n\beta_1^2N_2-n\beta_1^2+N_2-n)^2}
$,

$p_{\gamma_1,\gamma_1}^{\beta_1}
= \frac{(n+1)(\gamma_1^2n\beta_1+1)}{(\gamma_1^2n+1)^2},
$
\quad
$p_{\gamma_2,\gamma_2}^{\beta_1}
= \frac{n(n+1)(n\gamma_1^2+\beta_1)\gamma_1^2}{(n\gamma_1^2+1)^2},
$

$p_{\gamma_1,\gamma_2}^{\beta_1}
=p_{\gamma_2,\gamma_1}^{\beta_1}
= \frac{n(n+1)\gamma_1^2(1-\beta_1)}
{(n\gamma_1^2+1)^2}.
$\\
Proposition \ref{pro:3-1} (1) with
$\lambda=2$, $q=2$, $(l,k,j)=(0,0,0),
(0,0,1), (1,0,0), (1,0,1)$,
$(1,1,0),$ and 
$(1,1,1)$ imply

$p_{\beta_1,\beta_1}^{\beta_2}
=\frac{ N_2(1-\beta_1)(N_2-n-1)(N_2+n\beta_1-n)}{(2n\beta_1+n\beta_1^2N_2-\beta_1^2n+N_2-n)^2},
$

$p_{\beta_2,\beta_2}^{\beta_2}
= \frac{(1+N_2\beta_1-\beta_1)\bigg( n^2(N_2-1)(N_2-2)\beta_1^3+3n^2(N_2-2)\beta_1^2-3n(N_2-2n)\beta_1-(N_2-2n)(N_2-n)
\bigg)}{(2n\beta_1+n\beta_1^2N_2-\beta_1^2n+N_2-n)^2},
$

$p_{\beta_2,\beta_1}^{\beta_2}=
p_{\beta_1,\beta_2}^{\beta_2}= \frac{\beta_1N_2^2(n\beta_1+1)(N_2-n-1)}{(2n\beta_1+n\beta_1^2N_2-\beta_1^2n+N_2-n)^2},
$

$p_{\gamma_1,\gamma_1}^{\beta_2}
= -\frac{(n+1)(\gamma_1^2n\beta_1-\gamma_1^2n-N_2\beta_1+N_2\gamma_1^2+\beta_1-1)}{(\gamma_1^2n+1)^2(1+N_2\beta_1-\beta_1)},
$

$p_{\gamma_2,\gamma_2}^{\beta_2}
= \frac{(n+1)(\gamma_1^2n^2N_2\beta_1-\gamma_1^2n^2\beta_1+\gamma_1^2n^2+n-n\beta_1-N_2)\gamma_1^2}{(\gamma_1^2n+1)^2(1+N_2\beta_1-\beta_1)},
$

$p_{\gamma_1,\gamma_2}^{\beta_2}
=p_{\gamma_2,\gamma_1}^{\beta_2}
=\frac{ \gamma_1^2N_2(n+1)(n\beta_1+1)}{(\gamma_1^2n+1)^2(1+N_2\beta_1-\beta_1)},
$\\
Proposition \ref{pro:3-1} (2) with
$(\lambda,\mu)=(1,2)$,  $q=1$, $(l,k,j)=(0,0,0),
(0,0,1), (1,0,0), (0,1,0), (0,1,1)$, and
$(1,1,0)$ imply

$p_{\alpha_1,\gamma_1}^{\gamma_1}
= \frac{n(1-\gamma_1^2)}{\gamma_1^2n+1},
$
\quad
$p_{\gamma_1,\beta_1}^{\gamma_1}
= \frac{N_2(N_2-n-1)(\gamma_1^2n\beta_1+1)}{(\gamma_1^2n+1)(2n\beta_1+n\beta_1^2N_2-\beta_1^2n+N_2-n)},
$

$p_{\gamma_1,\beta_2}^{\gamma_1}
= -\frac{n(1+N_2\beta_1-\beta_1)(\gamma_1^2n\beta_1-\gamma_1^2n-N_2\beta_1+N_2\gamma_1^2+\beta_1-1)}{(\gamma_1^2n+1)(2n\beta_1+n\beta_1^2N_2-\beta_1^2n+N_2-n)},
$

$p_{\gamma_2,\beta_2}^{\gamma_1}
=\frac{ (1+N_2\beta_1-\beta_1)(n\beta_1+1)N_2\gamma_1^2n}{(\gamma_1^2n+1)(2n\beta_1+n\beta_1^2N_2-\beta_1^2n+N_2-n)},
$

$p_{\alpha_1,\gamma_2}^{\gamma_1}
= \frac{n(n+1)\gamma_1^2}{n\gamma_1^2+1},
$
\quad
$p_{\gamma_2,\beta_1}^{\gamma_1}
=\frac{\gamma_1^2nN_2(1-\beta_1)(N_2-n-1)}{(\gamma_1^2n+1)(2n\beta_1+n\beta_1^2N_2-\beta_1^2n+N_2-n)},
$\\
Proposition \ref{pro:3-1} (2) with
$(\lambda,\mu)=(1,2)$,  $q=2$, $(l,k,j)=(0,0,0),
(0,0,1), (1,0,0), (0,1,0), (0,1,1)$, and
$(1,1,0)$ imply

$p_{\alpha_1,\gamma_1}^{\gamma_2}
=\frac{n+1}{\gamma_1^2n+1}(=p_{\gamma_1,\gamma_1}^{\beta_0}),
$

$p_{\gamma_1,\beta_1}^{\gamma_2}
=\frac{ N_2(1-\beta_1)(N_2-n-1)}{(\gamma_1^2n+1)(2n\beta_1+n\beta_1^2N_2-\beta_1^2n+N_2-n)},
$

$p_{\gamma_2,\beta_2}^{\gamma_2}
= \frac{(1+N_2\beta_1-\beta_1)(\gamma_1^2n^2N_2\beta_1-\gamma_1^2n^2\beta_1+\gamma_1^2n^2+n-n\beta_1-N_2)}{(\gamma_1^2n+1)(2n\beta_1+n\beta_1^2N_2-\beta_1^2n+N_2-n)},
$

$p_{\alpha_1,\gamma_2}^{\gamma_2}
=\frac{ (n\gamma_1-1)(n\gamma_1+1)}{(\gamma_1^2n+1)},
$

$p_{\gamma_2,\beta_1}^{\gamma_2}
= \frac{N_2(N_2-n-1)(\gamma_1^2n+\beta_1)}{(\gamma_1^2n+1)(2n\beta_1+n\beta_1^2N_2-\beta_1^2n+N_2-n)},
$

$p_{\gamma_1,\beta_2}^{\gamma_2}
= \frac{(n\beta_1+1)N_2(1+N_2\beta_1-\beta_1)}{(\gamma_1^2n+1)(2n\beta_1+n\beta_1^2N_2-\beta_1^2n+N_2-n)}.
$\\

Thus we obtained the intersection numbers 
interms of $n$, $N_2$, $\beta_1$ and $\gamma_1$.
From the remaining equations given in 
Proposition \ref{pro:3-1} we obtained 
the following seven equalities between $n$, $N_1$, $\beta_1$,
$\gamma_1$, $w$ and $r$.
More precisely, 
Proposition \ref{pro:3-1} (1) with
$\lambda=1$, $q=0$, $(l,k,j)=(1,2,0)$ and $(3,1,0)$
imply the following (\ref{equ:n+1(1)}) and
 (\ref{equ:n+1(2)}) respectively. 
  Proposition \ref{pro:3-1} (1) with
$\lambda=2$, $q=0$, $(l,k,j)=(2,1,0)$
implies the following (\ref{equ:n+1(3)}).
\begin{eqnarray}
&&N_2(n\gamma_1^2-1)r^3w_2
+(n^2-1)\gamma_1=0,
\label{equ:n+1(1)}
\\%(1)
&&N_2\bigg(n^2(n+2)\gamma_1^4
-2n(2n+1)\gamma_1^2+n+2\bigg)r_2^4w_2\nonumber\\
&&\qquad+(n-1)(n-2)(n+1)^2\gamma_1^2=0,
\label{equ:n+1(2)}
\\%(2)
&& N_2\bigg(n(N_2-n-1)\beta_1^2+n(n-1)\beta_1+2n-N_2\bigg)r_2^3w_2\gamma_1
\nonumber\\
&&
\qquad+(n+1)(\gamma_1^2n-1)(1+N_2\beta_1-\beta_1)=0,
\label{equ:n+1(3)}%(3)
 \end{eqnarray}
 Since $n\geq 2$ and $\gamma_1>0$,
 (\ref{equ:n+1(1)}) implies
 $\gamma_1\not=\frac{1}{\sqrt{n}}$.
 Then  (\ref{equ:n+1(1)}), (\ref{equ:n+1(3)})
 and  (\ref{equ:n+1-beta2})
 imply 
\begin{eqnarray}&& 
\beta_1=\frac{1}{2\gamma_1^2n(n-1)(N_2-n-1)}
\times\nonumber\\
&&
\bigg[n^2(N_2-1)\gamma_1^4-n(2N_2+n^2-2n-1)\gamma_1^2+N_2-1
\nonumber\\
&&
+\bigg\{n^4(N_2-1)^2\gamma_1^8
-2n^3\bigg(2N_2^2+(n^2-4n-1)N_2+n^2+2n-1\bigg)\gamma_1^6\nonumber\\
&&
+n\bigg(2(2n^2-n+2)N_2^2
-(8n^3-4n^2+4n+4)N_2
+n^5+4n^4+2n^3-4n^2+3n
\bigg)\gamma_1^4\nonumber\\
&&
-2n\bigg(2N_2^2+(n^2-4n-1)N_2+n^2+2n-1\bigg)\gamma_1^2
+(N_2-1)^2\bigg\}^{\frac{1}{2}}
\bigg],\label{equ:n+1-beta1}\\
&& 
\beta_2=\frac{1}{2\gamma_1^2n(n-1)(N_2-n-1)}
\times\nonumber\\
&&
\bigg[n^2(N_2-1)\gamma_1^4-n(2N_2+n^2-2n-1)\gamma_1^2+N_2-1
\nonumber\\
&&
-\bigg\{n^4(N_2-1)^2\gamma_1^8
-2n^3\bigg(2N_2^2+(n^2-4n-1)N_2+n^2+2n-1\bigg)\gamma_1^6\nonumber\\
&&
+n\bigg(2(2n^2-n+2)N_2^2
-(8n^3-4n^2+4n+4)N_2
+n^5+4n^4+2n^3-4n^2+3n
\bigg)\gamma_1^4\nonumber\\
&&
-2n\bigg(2N_2^2+(n^2-4n-1)N_2+n^2+2n-1\bigg)\gamma_1^2
+(N_2-1)^2\bigg\}^{\frac{1}{2}}
\bigg].\label{equ:n+1-beta2-2}
\end{eqnarray}
Since $X_2$ is a strongly regular graph,
ratio of the squares of usual Euclidean distances between the 
points in $X_2$ is given by $\frac{k-1}{k}$ 
and $\left(\frac{2-\beta_1-\beta_2}
{\beta_1-\beta_2}\right)=2k-1$, where $k$ is an integer
satisfying
$k\geq 2$ (see \cite{B-B-2,B-B-3,L-R-S}). We can express 
$\left(\frac{2-\beta_1-\beta_2}
{\beta_1-\beta_2}\right)^2$ as follows.
\begin{eqnarray}
&&\left(\frac{2-\beta_1-\beta_2}
{\beta_1-\beta_2}\right)^2=
\bigg(n^2(N_2-1)\gamma_1^4-n(2nN_2-n^2-2n+1)\gamma_1^2+N_2-1\bigg)^2\times
\nonumber\\
&&
\bigg\{n^4(N_2-1)^2\gamma_1^8
-2n^3(2N_2^2+(n^2-4n-1)N_2+n^2+2n-1)\gamma_1^6
\nonumber\\
&&
+n\bigg((4n^2-2n+4)N_2^2-4(2n+1)(n^2-n+1)N_2
+n(n^4+4n^3+2n^2-4n+3)
\bigg)\gamma_1^4
\nonumber\\
&&
-2n\bigg(2N_2^2+(n^2-4n-1)N_2+n^2+2n-1\bigg)\gamma_1^2+(N_2-1)^2\bigg\}^{-1}
\label{equ:n+1-ratio}
\end{eqnarray}
Let $m=p_{\gamma_1,\gamma_1}^{\beta_0}$.
Then (\ref{equ:n+1-m1}) implies
$\gamma_1= \sqrt{\frac{n+1-m}{mn}}$. 
Since $\gamma_1>0$ and $s_{1,2}=2$, we must have 
$1\leq m\leq n$.
Then (\ref{equ:n+1(1)}) and (\ref{equ:n+1(2)})
imply 
\begin{eqnarray}
&&r_2=\frac{(n-2)(2m-n-1)\sqrt{m}\sqrt{n+1-m}}
{(-n^2-3n+6nm-6m^2-2+6m)\sqrt{n}},\label{equ:n+1-r}\\
&&w_2=\frac{(-n^2-3n+6nm-6m^2-2+6m)^3(n+1)n(n-1)}{(n-2)^3N_2m(n+1-m)(n+1-2m)^4},\label{equ:n+1-w}
\end{eqnarray}
Since $w_2,\ r_2>0$, we must have
$$-n^2-3n+6nm-6m^2-2+6m>0, \qquad
2m-n-1>0.$$
Therefore 
$$\frac{n+1}{2}<m<\frac{n+1}{2}+\frac{\sqrt{3(n^2-1)}}{6}$$
holds. Then (\ref{equ:n+1-beta1}), (\ref{equ:n+1-beta2-2}),
(\ref{equ:n+1-ratio}), (\ref{equ:n+1-r}) and (\ref{equ:n+1-w})
imply
\begin{equation}
\left(\frac{2-\beta_1-\beta_2}
{\beta_1-\beta_2}\right)^2=
\frac{P_1(n,N_2,m)}{P_2(n,N_2,m)}\label{equ:4-13},
\end{equation}
where
\begin{eqnarray}
&&P_1(n,N_2,m)=n\bigg((2N_2-n-1)m^2+(n+1)(n+1-2N_2)m+(N_2-1)(n+1)\bigg)^2,\nonumber\\
&&\\
&&
P_2(n,N_2,m)=
\bigg(4N_2^2-4(n+1)N_2+n(n+1)^2)(m-2(n+1))m^3
\nonumber\\
&&
+\bigg(4(n^2+4n+1)N_2^2
-(2n^3+20n^2+22n+4)N_2
+n(n^2+2n+3)(n+1)^2
\bigg)m^2
\nonumber\\
&&
-2n(n+1)\bigg(4N_2^2+(n+1)(n-5)N_2
+(n+1)^2\bigg)m+n(n+1)^2(N_2-1)^2.
\end{eqnarray}
 Let $F(n,x,y)=\frac{P_1(n,x,y)}{P_2(n,x,y)}$
 and consider the behavior of $F(n,x,y)$
 for $\frac{n(n+1)}{2}+1\leq x\leq \frac{n(n+3)}{2}$,
 $\frac{n+1}{2}<y<\frac{n+1}{2}+\frac{\sqrt{3(n^2-1)}}{6} $.
 We have
 \begin{eqnarray}
&& \frac{\partial F(n,x,y)}{\partial y}
 =\frac{4nx(n-1)(x-n-1)(n+1-2y)}{P_2(n,x,y)^2}\times\nonumber\\
 &&
\bigg((2x-n-1)y^2+(n+1)(n+1-2x)y+(x-1)(n+1)\bigg)\times\nonumber\\
 &&
\bigg(-(n^2-n-2+2x)y^2+(n+1)(n^2-n-2+2x)y-n(x-1)(n+1)\bigg)
\nonumber\\
 &&\label{equ:dFy}
 \end{eqnarray}
  \begin{eqnarray}
&& \frac{\partial F(n,x,y)}{\partial x}
 =-\frac{4ny(y-1)(n^2-1)(n+1-y)(n-y)
}{P_2(n,x,y)^2}\times
 \nonumber\\
 &&\bigg(x(n(2y-1)-2y^2+2y-1)-(n+1)( y(n-y)+y-1)\bigg)\times
 \nonumber\\
 &&
 \bigg((n-1)x-(n+1)(y(n-y)+y-1)\bigg)
\label{equ:dFN2}
 \end{eqnarray}
 (\ref{equ:dFy}) and (\ref{equ:dFN2}) imply
$ \frac{\partial F(n,x,y)}{\partial y}<0$ and 
$ \frac{\partial F(n,x,y)}{\partial x}>0$
for any $x,y$ satisfying $ \frac{n(n+1)}{2}+1\leq x\leq 
\frac{n(n+3)}{2}$, $\frac{n+1}{2}< y<\frac{n+1}{2}+\frac{\sqrt{3(n^2-1)}}{6}$.
We can check the followings easily
\begin{eqnarray}&&F(n,x,y)<F\left(n,x,\frac{n+1}{2}\right)\leq
F\left(n,\frac{n(n+3)}{2},\frac{n+1}{2}\right)=n+3,\nonumber\\
&&F(n,x,y)>F\left(n,x,\frac{n+1}{2}+\frac{\sqrt{3(n^2-1)}}{6}\right)\nonumber\\
&&
\qquad> 
F\left(n,\frac{n(n+1)}{2},\frac{n+1}{2}
+\frac{\sqrt{3(n^2-1)}}{6}\right)
=n+2.
\nonumber
\end{eqnarray}
Therefore $F(n,x,y)$ cannot be an integer for any 
$x,y$ satisfying 
$ \frac{n(n+1)}{2}+1\leq x\leq 
\frac{n(n+3)}{2}$, $\frac{n+1}{2}< y<\frac{n+1}{2}+\frac{\sqrt{3(n^2-1)}}{6}$. Hence we must have
$N_2=\frac{n(n+1)}{2}$ and $(X, w)$ is a Euclidean tight $4$-design. 
This completes the proof of Theorem \ref{theo:4-1}.
\hfill\qed\\

In the following we assume $N_2\geq N_1\geq n+2$.
Hence we must have $s_{1,1}=s_{2,2}=2$.
Let $\alpha_i=\alpha_{1,1}^{(i)}$ for $i=0,1,2$
and assume $\alpha_1>\alpha_2$.
We will prove the following theorem.

\begin{theo}
\label{theo:4-2} 
Let definition and notation be given as above.
We have the following assertions.
\begin{enumerate}
\item The following hold.%(1)
\begin{enumerate}
\item
$\gamma_1\gamma_2=-\frac{1}{n}$.
\item $(N_1-1)\alpha_1+1\neq 0$
and 
\begin{equation}
\alpha_2=-\frac{n\alpha_1-n+N_1}
{n((N_1-1)\alpha_1+1)}.
\label{equ:alpha2}
\end{equation}
\item $(N_2-1)\beta_1+1\neq 0$ and 
\begin{equation}\beta_2 = -\frac{n\beta_1-n+N_2}{n((N_2-1)\beta_1+1)}.\label{equ:beta2}
\end{equation}
\item $p_{\gamma_1,\gamma_1}^{\alpha_0}
=\frac{N_2}{1+n\gamma_1^2}$.
\item $p_{\gamma_1,\gamma_1}^{\beta_0}
=\frac{N_1}{1+n\gamma_1^2}$.
\end{enumerate}
\item%(2)
 If $\gamma_1=\frac{1}{\sqrt{n}}$, then 
$(X,w)$ is similar to the 
Euclidean $4$-design given in Theorem \ref{theo:1-4}
(2) (i) or to the one given in Theorem \ref{theo:1-4} (2) (ii). 
\item%(3)
If $\gamma_1\not=\frac{1}{\sqrt{n}}$, then 
$(X,w)$ is a Euclidean tight $4$-design. 
\end{enumerate}
\end{theo}
Theorem \ref{theo:4-2} together with Theorem \ref{theo:4-1}
implies Theorem \ref{theo:1-4} (2). In the following 
we give the proof for Theorem \ref{theo:4-2}.\\

\noindent{\bf Proof for Theorem \ref{theo:4-2} (1) }\\
{\bf (i)}
The equations in Proposition \ref{pro:3-1}
for $\boldsymbol x=\boldsymbol y\in X_1$
imply $\gamma_1\gamma_2=-\frac{1}{n}$. \\
{\bf (ii)}
The equations in Proposition \ref{pro:3-1}
for $\boldsymbol x=\boldsymbol y\in X_1$ also
imply
$$-n((N_1-1)\alpha_1+1)\alpha_2-n\alpha_1+n-N_1=0.$$
If $(N_1-1)\alpha_1+1=0$, then we must have 
$-n\alpha_1+n-N_1=0$. This implies 
$N_1=n+1$ and contradicts the assumption
that $N_1\geq n+2$.
Hence we obtain  (ii).\\ 
{\bf (iii)}
The equations in Proposition \ref{pro:3-1}
for
$\boldsymbol x=\boldsymbol y\in X_2$,
implies
$$(n(N_2-1)\beta_1+n)\beta_2+n\beta_1+N_2-n=0.
$$
If $n(N_2-1)\beta_1+n=0$, then 
we must have $n\beta_1+N_2-n=0$. This implies
$N_2=n+1$. This is a contradiction.
Hence we have (iii).\\
{\bf (iv) and (v) }
We obtain (iv) and (v) using 
using the equations given in Proposition \ref{pro:3-1}.
Explicit formulas for the intersection numbers
are given in terms of $n, N_1, N_2, \alpha_1,
\beta_1, \gamma_1$.
The reader can find them in Appendix III.\\

Using the intersection numbers expressed
in terms of $n, N_1, N_2, \alpha_1,
\beta_1, \gamma_1$, we obtain the nine 
equations given below.
If $N_1, N_2, \alpha_1, \beta_1,
\gamma_1, w_2$ and $r_2$ satisfy all of the
nine equations
and if the intersection numbers satisfy
integral condition, then 
$(X,w)$ satisfies the conditions of
Euclidean $4$-design.
That is, we obtain feasible parameters for a Euclidean
$4$-design and the corresponding coherent configuration.

Proposition \ref{pro:3-1} (1) with
$(\lambda,q,l,k,j)=(1,0,3,0,0),(1,0,2,2,0),
(2,0,2,1,0),\\ (2,0,3,1,0),
(1,1,2,2,0),(1,2,2,2,0),
(2,1,2,2,0),(2,2,2,2,0)$
and 
and Proposition \ref{pro:3-1} (2) with
$(\lambda,\mu,q,l,k,j)=(1,2,1,2,2,0)$
imply the following
(\ref{equ:I}),
(\ref{equ:II}),
(\ref{equ:III}),
(\ref{equ:IV}),
(\ref{equ:V}),
(\ref{equ:VI}),
(\ref{equ:VII}),
(\ref{equ:VIII}) and
(\ref{equ:IX}) respectively.
\begin{eqnarray}&&N_2(N_1\alpha_1-\alpha_1+1)(n\gamma_1^2-1)w_2r_2^3\nonumber\\
&&\qquad
+\gamma_1N_1\bigg(n(N_1-n-1)\alpha_1^2
+n(n-1)\alpha_1+2n-N_1
\bigg)=0,\label{equ:I}\\%I, Test01(3,0,0)
&&N_2(N_1\alpha_1-\alpha_1+1)^2
(n^2\gamma_1^4(n+2)-2\gamma_1^2n(2n+1)+n+2)w_2r_2^4\nonumber\\
&&\qquad+\gamma_1^2N_1\bigg\{n^2(n+2)
(N_1-1)(N_1-n-1)\alpha_1^4\nonumber\\
&&\qquad+n\bigg(2(2n+1)N_1^2
-(n^3+4n^2+11n+2)N_1+3n(n^2+2n+1)\bigg)\alpha_1^2\nonumber\\
&&\qquad+2n(n-1)(n^2+n-2N_1)
\alpha_1+3n^2(n+1) +(n+2)N_1^2 -3n(n+2)N_1\bigg\}=0,\nonumber
\\
&&\label{equ:II}%II,Test01(2,2,0)
\end{eqnarray}
\begin{eqnarray}
&&\gamma_1N_2\bigg(
n(N_2-n-1)\beta_1^2
+n(n-1)\beta_1+2n-N_2\bigg)w_2r_2^3\nonumber\\
&&\qquad+N_1(n\gamma_1^2-1)
\bigg((N_2-1)\beta_1+1\bigg)=0,
\label{equ:III}%III, Test02(2,1,0)
\end{eqnarray}
\begin{eqnarray}
&&\gamma_1^2N_2
\bigg\{n^2(n+2)(N_2-1)(N_2-n-1)\beta_1^4\nonumber\\
&&\qquad
-n\bigg(2(2n+1)N_2^2-(2+4n^2+11n+n^3)N_2+3n(n+1)^2\bigg)\beta_1^2
\nonumber\\
&&\qquad
+2n(n-1)(n^2+n-2N_2)\beta_1
+(n+2)N_2^2-3n(n+2)N_2+3n^2(n+1)\bigg\}w_2r_2^4\nonumber\\
&&\qquad+N_1(N_2\beta_1+1-\beta_1)^2
\bigg(n+2-2n(2n+1)\gamma_1^2+(n+2)n^2\gamma_1^4\bigg)=0,\label{equ:IV}%IV Test02(3,1,0)
\end{eqnarray}
\begin{eqnarray}
&&-N_2(N_1\alpha_1-\alpha_1+1)^2\bigg(
n^2\alpha_1(n+2)\gamma_1^4
-2n(n\alpha_1^2+(n+2)\alpha_1-1)\gamma_1^2+\alpha_1(n+2)
\bigg)w_2r_2^4\nonumber\\
&&+N_1\gamma_1^2\bigg(
n^2(N_1-1)(n+2)(2n-N_1+1)\alpha_1^5
\nonumber\\
&&-n^2(2nN_1+2+2n^2N_1-4n-3n^2-2N_1^2)\alpha_1^4
+2n(N_1-1)((n+2)N_1-3n^2-4n)\alpha_1^3\nonumber\\
&&-n(-2n^2N_1+2N_1^2+n-2nN_1+4n^2+n^3)\alpha_1^2
\nonumber\\
&&+(-n^3-4n^2+4n(n+1)N_1-(n+2)N_1^2)\alpha_1+n^2
\bigg)=0,\label{equ:V}%V Testq1(2,2,0,1)
\end{eqnarray}
\begin{eqnarray}
&&N_2\bigg\{n\bigg(
-n^2(N_1-1)(n+2)\gamma_1^4+2n(-2n+N_1^2+N_1n-1)\gamma_1^2-(N_1-1)(n+2)
\bigg)\alpha_1^2\nonumber\\
&&+\bigg(
n^2(n+2)(-2n+N_1n-N_1^2+N_1)\gamma_1^4\nonumber\\
&&
-2n(-2n+2N_1-4n^2+3N_1n-N_1^2n-2N_1^2+N_1n^2)\gamma_1^2\nonumber\\
&&
+(n+2)(-2n+N_1n-N_1^2+N_1)
\bigg)\alpha_1
+n^2(n+2)(n-N_1)\gamma_1^4\nonumber\\
&&
-2n(-2N_1+2n^2+n-3N_1n+N_1^2)\gamma_1^2+(n+2)(n-N_1)
\bigg\}
(N_1\alpha_1-\alpha_1+1)w_2r_2^4\nonumber\\
&&
+\bigg(
n^3(N_1-1)(n+2)(n-N_1+1)\alpha_1^5\nonumber\\
&&
-n^2(n+2)(-n+2N_1n-N_1^2+N_1)(n-N_1+1)\alpha_1^4\nonumber\\
&&
+n^2(2N_1^3-4N_1+6n^2+3n+2N_1^2+3n^3+N_1n^3-11N_1n-2N_1^2n^2
)\alpha_1^3\nonumber\\
&&
-n(-26N_1n^2+16N_1^2n-4N_1^3+4N_1^2-4N_1n+6n^3+n^2+5n^4+8N_1^2n^2-8N_1n^3-2N_1^3n)\alpha_1^2\nonumber\\
&&
+n(12N_1n-4N_1^2-3n^3+5N_1^2n-3N_1n^2-2N_1^3-5n^2+2n^4+2N_1^2n^2-4N_1n^3)\alpha_1\nonumber\\
&&
-7N_1n^3+8N_1^2n+3n^3+5N_1^2n^2-2N_1^3-N_1^3n-10N_1n^2+3n^4\bigg)N_1\gamma_1^2=0,
\label{equ:VI}% VI Testq1(2,2,0,2)
\end{eqnarray}
\begin{eqnarray}
&&N_2\bigg(n^2\beta_1^5(N_2-1)(2+n)(-2n+N_2-1)-n^2(2N_2^2-2n^2N_2-2N_2n+3n^2+4n-2)\beta_1^4
\nonumber\\
&&
-2n(N_2-1)(2N_2+N_2n-3n^2-4n)\beta_1^3
+n(-2N_2n+2N_2^2-2n^2N_2+n+n^3+4n^2)\beta_1^2
\nonumber\\
&&
+(-4n^2N_2+4n^2-4N_2n+nN_2^2+n^3+2N_2^2)\beta_1-n^2\bigg)\gamma_1^2w_2r_2^4
\nonumber\\
&&
+N_1\bigg(n^2(n+2)\beta_1\gamma_1^4
-2n(n\beta_1^2+(n+2)\beta_1-1)\gamma_1^2+(n+2)\beta_1\bigg)\bigg((N_2-1)\beta_1+1\bigg)^2=0,
\nonumber\\
&&\label{equ:VII}%VII Testq2(2,2,0,1)
\end{eqnarray}
%\newpage
\begin{eqnarray}
&&N_2\gamma_1^2
\bigg(-n^3(N_2-1)(n+2)(N_2-n-1)\beta_1^5
\nonumber\\
&&-n^2(n+2)(n-2N_2n+N_2^2-N_2)(N_2-n-1)\beta_1^4
\nonumber\\
&&+n^2(3n^3+n^3N_2-11N_2n+2N_2^2+3n-2n^2N_2^2+2N_2^3+6n^2-4N_2)\beta_1^3
\nonumber\\
&&+n(2(n+2)N_2^3
-4(2n^2+4n+1)N_2^2
+n(8n^2+26n+4)N_2-(5n^2+6n+1)n^2)\beta_1^2
\nonumber\\
&&-n(-2n^4+3n^2N_2+4N_2^2-2n^2N_2^2+5n^2-5nN_2^2-12N_2n+3n^3+4n^3N_2+2N_2^3)\beta_1
\nonumber\\
&&-N_2^3n-7n^3N_2+5n^2N_2^2-2N_2^3-10n^2N_2+3n^4+3n^3+8nN_2^2\bigg)w_2r_2^4
\nonumber\\
&&
+N_1\bigg\{-n\bigg(n^2(N_2-1)(2+n)\gamma_1^4-2n(-1+N_2n-2n+N_2^2)\gamma_1^2
\nonumber\\
&&
+(N_2-1)(2+n)\bigg)
(N_2-1)\beta_1^3
-\bigg(n^2(N_2-1)(2+n)(-N_2n+3n+N_2^2-N_2)\gamma_1^4
\nonumber\\
&&
-2n(-3n-6n^2-n^2N_2^2+2N_2^3+N_2^3n+6n^2N_2-3nN_2^2+5N_2n+2N_2-4N_2^2)\gamma_1^2
\nonumber\\
&&
+(N_2-1)(2+n)(-N_2n+3n+N_2^2-N_2)\bigg)\beta_1^2
\nonumber\\
&&
-\bigg(n^2(2+n)(3n-2N_2n+2N_2^2-2N_2)\gamma_1^4
\nonumber\\
&&
+2n(-6n^2-3n+N_2^3-5N_2^2+7N_2n+4N_2+3n^2N_2-4nN_2^2)\gamma_1^2
\nonumber\\
&&
+(2+n)(3n-2N_2n+2N_2^2-2N_2)\bigg)\beta_1
\nonumber\\
&&
-n^2(2+n)(N_2-n)\gamma_1^4-2n(N_2^2-3N_2n+2n^2+n-2N_2)\gamma_1^2-(2+n)(N_2-n)\bigg\}=0,
\nonumber\\
&&\label{equ:VIII}%VIII Testq2(2,2,0,2)
\\
&&(n\gamma_1^2-1)
\bigg\{(N_1\alpha_1-\alpha_1+1)\bigg(-n(n+2)(N_2-n-1)\beta_1^2+n(-n^2-n+2N_2)\beta_1
\nonumber\\
&&
+N_2n-2n^2-2n+2N_2\bigg)N_2r_2^4w_2
\nonumber\\
&&+(1+N_2\beta_1-\beta_1)\bigg(-n(n+2)(N_1-n-1)\alpha_1^2+n(-n^2-n+2N_1)\alpha_1
\nonumber\\
&&
+N_1n-2n^2-2n+2N_1\bigg)N_1
\bigg\}=0,\label{equ:IX}%IX TestDq(2,2,0,2)
\end{eqnarray}

\noindent
{\bf Proof for Theorem \ref{theo:4-2} (2)}\\
Since $\gamma_1=\frac{1}{\sqrt{n}}$ (\ref{equ:I}) and (\ref{equ:alpha2}) imply
$$\alpha_1=\frac{-n^2+n
+\sqrt{n\bigg(4N_1^2-4(3n+1)N_1+n^3+6n^2+9n\bigg)}}
{2n(N_1-n-1)},$$ 
$$\alpha_2=\frac{-n^2+n
-\sqrt{n\bigg(4N_1^2-4(3n+1)N_1+n^3+6n^2+9n\bigg)}}
{2n(N_1-n-1)}.$$
Also  (\ref{equ:III}) and (\ref{equ:beta2}) imply
$$\beta_1= \frac{-n^2+n
+\sqrt{n\bigg(4N_2^2-4(3n+1)N_2+n^3+6n^2+9n\bigg)}}
{2n(N_2-n-1)},
$$
$$\beta_2= \frac{-n^2+n
-\sqrt{n\bigg(4N_2^2-4(3n+1)N_2+n^3+6n^2+9n
\bigg)}}{2n(N_2-n-1)}
$$
Hence we obtain
\begin{eqnarray}
&&\left(\frac{2-\alpha_1-\alpha_2}{\alpha_1-\alpha_2}\right)^2
=\frac{(n-2N_1+3)^2n}{(9n+6n^2+n^3-12N_1n+4N_1^2-4N_1)}
\nonumber\\
&&
\left(\frac{2-\beta_1-\beta_2}{\beta_1-\beta_2}\right)^2
=\frac{(-2N_2+n+3)^2n}{9n+6n^2+n^3-4N_2-12nN_2+4N_2^2}
\nonumber
\end{eqnarray}
Since $N_1+N_2\geq \frac{(n+2)(n+1)}{2}$,
$N_2\geq N_1\geq n+2$, we must have $N_2\geq \frac{(n+2)(n+1)}{4}$.
Let $F(n,x)=\frac{(-2x+n+3)^2n}{9n+6n^2+n^3-4x-12nN_2+4x^2}$
and consider the behavior of $F(n,x)$.
$$\frac{d F(n,x)}{dx}=
\frac{4(2x-n-3)n(n-1)(n^2+4n+3-4x)}
{(9n+6n^2+n^3-4x-12nx+4x^2)^2
}$$
For $x\in [\frac{(n+2)(n+1)}{4}, \frac{n(n+3)}{2}]$
$F(x)$ takes the maximal value at
$x=\frac{n^2+4n+3}{4}$.
Moreover we have
 $$F\left(n,\frac{n^2+4n+3}{4}\right)=\frac{n(n+3)}{n-1}=
 n+4+\frac{4}{n-1},$$
 $$F\left(n,\frac{(n+2)(n+1)}{4}\right)=
 n+4+
 \frac{4n^3-8n^2-28n+16}{12n-3n^2-2n^3-4+n^4}
 $$
 $$F\left(n,\frac{n(n+3)}{2}\right)=
 n+3.
 $$
Since $\frac{(n+2)(n+1)}{4}\leq N_2
\leq \frac{n(n+3)}{2}$, we have the following 
\begin{eqnarray}
&&\left(\frac{2-\beta_1-\beta_2}{\beta_1-\beta_2}\right)^2
=n+3,\ \mbox{or}\ n+4\ 
\mbox{for $n\geq 6$,}
\end{eqnarray}
where $n+3$ (or $n+4$ respectively) is
the square of an odd integer.
\begin{eqnarray}
&&\left(\frac{2-\beta_1-\beta_2}{\beta_1-\beta_2}\right)^2=9
\ \mbox{for $n\leq 5$}.
\end{eqnarray}

\noindent
If $n\geq 6$ and $(\frac{2-\beta_1-\beta_2}{\beta_1-\beta_2})^2=n+3$
hold,
then we must have
$N_2=\frac{n(n+3)}{2}$.
Since $n+3=(2k-1)^2$ with an integer $k\geq 2$,
we have 
$N_2=(2k^2-2k-1)(2k-1)^2$ which is an odd integer.
However equations in Proposition \ref{pro:3-1} (1)
with $\lambda=1$ and $q=0$ implies
$p_{\gamma_1,\gamma_1}^{\alpha_0}=\frac{N_2}{2}$.
Hence $N_2$ must be an even integer.
This is a contradiction. 
Hence if $n\geq 6$, we only need to 
consider the case where $n+4$ is the
square of an odd integer
and 
$\left(\frac{2-\beta_1-\beta_2}{\beta_1-\beta_2}\right)^2=n+4$.
 
If $n\leq 5$, then we must have
$\left(\frac{2-\beta_1-\beta_2}{\beta_1-\beta_2}\right)^2=9$.
This implies $n=2$ or $n=5$.
If  $n=2$, then elementary computations imply 
that $(X,w)$ is similar to the one given in Theorem \ref{theo:1-4}
(2) (i).

If $n=5$, then $\left(\frac{2-\beta_1-\beta_2}{\beta_1-\beta_2}\right)^2=9(=n+4)$ holds.

In the following we may assume 
$n\geq 5$ and  
$(\frac{2-\beta_1-\beta_2}{\beta_1-\beta_2})^2=n+4=(2k-1)^2$.
Then $F(n,N_2)=n+4$ imples 
$N_2=\frac{n^2+5n+2}{4}+\frac{(n-1)\sqrt{n+4}}{4}
=2k^3(2k-3)$.
Then we must have 
$N_1\geq \frac{(n+2)(n+1)}{2}-N_2=
(2k+1)(2k^3-6k^2+4k+1)
$.
Let 
\begin{eqnarray}&&G(k,x)=F((2k-1)^2-4,x)\nonumber\\
&&\qquad=\frac{
(x-2k^2+2k)^2(2k+1)(2k-3)}{(16k^6-48k^5+36k^4+8k^3-12k^2-12xk^2+12xk+8x+x^2)}
\end{eqnarray}
Then 
\begin{eqnarray}&&\frac{\partial G(k,x)}{\partial x}=
\nonumber\\
&&\frac{8(x-2k^2+2k)(2k+1)(2k-3)(k^2-k-1)(4k^4-8k^3+2k^2+2k-x)
}{(16k^6-48k^5+36k^4+8k^3-12k^2-12xk^2+12xk+8x+x^2)^2
}
\end{eqnarray}
Then for $x\in [(2k+1)(2k^3-6k^2+4k+1),2k^3(2k-3)]
$
$G(k,x)$ takes the maximal value at
$x=4k^4-8k^3+2k^2+2k
(=\frac{n^2+4n+3}{4})$, $G(k,4k^4-8k^3+2k^2+2k)=
n+4+\frac{4}{n-1}
$,
$G(k,2k^3(2k-3))=(2k-1)^2(=n+4)$
and 
$$G(k,(2k+1)(2k^3-6k^2+4k+1))=
\frac{(4k^4-10k^3+8k+1)^2}
{4k^6-16k^5+8k^4+28k^3-12k^2-20k-3}.
$$

If $k\geq 4$, then $G(k,(2k+1)(2k^3-6k^2+4k+1))>n+3$.
Hence we must have $G(k,N_1)=n+4=(2k-1)^2$.

For $k=2$ and $3$, case by case computations
imply that we must also have $G(k,N_1)=n+4=(2k-1)^2$.

Then we must have 
 $ N_1= 2k^3(2k-3)(=N_2)$ or $2(2k+1)(k-1)^3$.

If $N_1=N_2= 2k^3(2k-3)$.
Then 
$w_2=\frac{k(2k-3)}{r_2^4(2k+1)(k-1)}$,
$\alpha_1=\beta_1=\frac{1}{2k+1}$,
$\alpha_2=\beta_2=-\frac{k+1}{(2k+1)(k-1)}$.
Then the second of the nine equation, (\ref{equ:II}),
implies 
$$-\frac{128}{2k+1}(k^2-k-1)(k-1)^3(2k-3)^2k^7
=0.$$
This is impossible.
Hence we must have $N_1=2(2k+1)(k-1)^3<N_2$.

$w_2=\frac{(2k+1)^2(k-1)^4}{(2k-3)^2k^4
}r_2^{-4}
$.

$\alpha_1=\frac{k-2}{k(2k-3)},\
\alpha_2=-\frac{1}{2k-3}
$,

$\beta_1=\frac{1}{2k+1},\
\beta_2=-\frac{k+1}{(2k+1)(k-1)},
$

$\gamma_1=\frac{1}{\sqrt{4k^2-4k-3}},\
\gamma_2=-\frac{1}{\sqrt{4k^2-4k-3}}.$\\
Thus $(X,w)$ is similar to the one having the parameter
given in Theorem \ref{theo:1-4} (2) (ii).
This completes the proof for Theorem \ref{theo:4-2} (2).\\

\noindent
{\bf Proof for Theorem \ref{theo:4-2} (3)}\\
Let $m_1=p_{\gamma_1,\gamma_1}^{\beta_0}$
and $m_2=p_{\gamma_1,\gamma_1}^{\alpha_0}$.
Then rquations in (iv) and (v) of Theorem \ref{theo:4-2} (1)
imply
 \begin{equation}\gamma_1=\sqrt{\frac{N_1-m_1}{nm_1}}
=\sqrt{\frac{N_2-m_2}{nm_2}},\ \mbox{and then}\
\frac{N_1}{N_2}=\frac{m_1}{m_2}.
\label{equ:gamma-m}
\end{equation}
Let $W=\frac{N_2}{N_1}r_2^3w_2$.
Then (\ref{equ:I}) and (\ref{equ:III}) 
imply
the following equations.
\begin{eqnarray}
&&(N_1\alpha_1-\alpha_1+1)(n\gamma_1^2-1)W
\nonumber\\
&&
\qquad\qquad
 +\gamma_1(n(N_1-n-1)\alpha_1^2
+n(n-1)\alpha_1-N_1+2n)=0
\end{eqnarray}
and
\begin{eqnarray}
&&(N_2\beta-\beta+1)(n\gamma_1^2-1)W^{-1}
\nonumber\\
&&\qquad\qquad
+\gamma_1(n(N_2-n-1)\beta_1^2+n(n-1)\beta_1-N_2+2n)
=0.
\end{eqnarray}
Then we obtain
\begin{eqnarray}
&&\alpha_1=\frac{(N_1-1)(1-n\gamma_1^2)W
-n(n-1)\gamma_1+\sqrt{D_{\alpha}}}
{2n(N_1-n-1)\gamma_1},\label{equ:alpha1}\\
&&\mbox{where}\quad
D_{\alpha}=
(N_1-1)^2(1-n\gamma_1^2)^2
W^2
-2n(1-n\gamma_1^2)(n+3+nN_1-3N_1)\gamma_1W
\nonumber\\
&&+n(4N_1^2-4(3n+1)N_1+n(n+3)^2)
\gamma_1^2\nonumber
\end{eqnarray}
and 
\begin{eqnarray}
&&\beta_1=\frac{(N_2-1)(1-n\gamma_1^2)W^{-1}
-n(n-1)\gamma_1+\sqrt{D_{\beta}}}
{2n(N_2-n-1)\gamma_1},\label{equ:beta1}\\
&&\mbox{where}\quad
D_{\beta}=
(N_2-1)^2(1-n\gamma_1^2)^2
W^{-2}
-2(1-n\gamma_1^2)(n+3+nN_2-3N_2)\gamma_1W^{-1}
\nonumber\\
&&+n(4N_2^2-4(3n+1)N_2+n(n+3)^2)
\gamma_1^2\nonumber
\end{eqnarray}

Then (\ref{equ:IX}), (\ref{equ:alpha1}) and
(\ref{equ:beta1}) imply
\begin{eqnarray}
&&\bigg(
n(nN_1-3N_1+n+3)\gamma_1
-(N_1-1)^2(1-n\gamma_1^2)W
-(N_1-1)\sqrt{D_{\alpha}}\bigg)
\times
\nonumber\\
&&
\bigg(n(nN_2-3N_2+n+3)\gamma_1W
-(N_2-1)^2(1-n\gamma_1^2)
-(N_2-1)\sqrt{D_{\beta}}\bigg)
\times
\nonumber\\
&&
\bigg(
2n\gamma_1-(n+2)(1-n\gamma_1^2)W
+r_2(2n\gamma_1W-(n+2)(1-n\gamma_1^2)))
\bigg)=0,
\nonumber\\
&&
\end{eqnarray}
where $D_{\alpha}$ and $D_{\beta}$ are given in
(\ref{equ:alpha1}) and (\ref{equ:beta1}) respectively.
Since $N_2\geq N_1\geq n+2$ and $\gamma_1>0$,
we must have
\begin{eqnarray}
&&2n\gamma_1-(n+2)(1-n\gamma_1^2)W
+r_2
\bigg(2n\gamma_1W
-(n+2)(1-n\gamma_1^2)\bigg)
=0.
\end{eqnarray}
If $2n\gamma_1W
-(n+2)(1-n\gamma_1^2)
=0$, then
we must have $1-n\gamma_1^2>0$
and
$2n\gamma_1-(n+2)(1-n\gamma_1^2)W
=0$.
Then we obtain
$$n^2(n+2)^2\gamma_1^4
-2n(n^2+6n+4)\gamma_1^2
+(n+2)^2=0.
$$
Since $\gamma_1^2$ must be a rational
number, $n(n+1)(n+4)$ must be a square 
of an positive integer.
If $n=2$, then $\gamma_1=\frac{1}{2}$
and $\gamma_2=-1$.
Then we obtain $W=1$ and
$\beta_1=1$. Hence
this case does not occur.
For $n\geq 3$, the following proposition shows that
this case does not occur.
\noindent
\begin{pro}\label{pro:Kaneko1}
For any integer $n\geq 3$, $n(n+1)(n+4)$
cannot be the square of an integer.
\end{pro}
{\bf Proof} Kaneko \cite{K}.\\

\noindent
Next, we assume $2n\gamma_1W-(n+2)
(1-n\gamma_1^2)\neq 0$.
Then we must have
\begin{eqnarray}
&&r_2=\frac{
(n+2)(1-n\gamma_1^2)W-2n\gamma_1}
{2n\gamma_1W-(n+2)
(1-n\gamma_1^2)}
\end{eqnarray}
Since $\gamma_1>0$ and $r_2>0$, we must have 
$1-n\gamma_1^2 > 0$, hence $0<\gamma_1 <\frac{1}{\sqrt{n}}$.

Then (\ref{equ:II}) implies
$$\left((n\gamma_1^2-1)W+(n-1)\gamma_1\right)
\left(P_0-(N_1-1)P_1\sqrt{D_{\alpha}}\right)P=0$$
where
\begin{eqnarray}
&&P=2(N_1-1)(1-n\gamma_1^2)\gamma_1W^2
-\bigg(n^2(n+2)\gamma_1^4
-6n\gamma_1^2+n+2\bigg)W
\nonumber\\
&&
+(n^2+3n-2N_1)(1-n\gamma_1^2)\gamma_1,
\nonumber\\
&&
P_0=
(N_1-1)^4(1-n\gamma_1^2)^2W^2
-2n(N_1-1)^2(n+3+nN_1-3N_1)(1-n\gamma_1^2)\gamma_1W
\nonumber\\
&&
+n\bigg(2N_1^4-6(n+1)N_1^3
+(n^3+6+21n)N_1^2
-(6n^2+24n+2)N_1
+n^3+6n^2+9n\bigg)\gamma_1^2
\nonumber\\
&&
\nonumber\\
&&
P_1=\bigg((N_1-1)^2(n\gamma_1^2-1)W
+n(nN_1-3N_1+n+3)\gamma_1\bigg),
\nonumber\\
\end{eqnarray}
 If $W=\frac{(n-1)\gamma_1}{1-n\gamma_1^2}$,
then (\ref{equ:alpha1})
implies $\alpha_1=-\frac{1}{n}$ and $\alpha_2=1$. 
This is a contradiction.
On the other hand we have
$P_0^2-(N_1-1)^2P_1^2D_{\alpha}=  4n^2N_1^2(N_1-n-1)^6W^4\gamma_1^4>0$.
Hence we must have
\begin{eqnarray}
&&2(N_1-1)(1-n\gamma_1^2)\gamma_1W^2
-\bigg((n+2)(n^2\gamma_1^4+1)
-6n\gamma_1^2\bigg)W
\nonumber\\
&&
+(n^2+3n-2N_1)(1-n\gamma_1^2)\gamma_1=0.
\label{equ:WP}
\end{eqnarray}
If we use (\ref{equ:IV}) instead of (\ref{equ:II}), then we obtain the following.
 $$((n\gamma_1^2-1)+(n-1)\gamma_1W)
 (Q_0-(N_2-1)Q_1\sqrt{D_{\beta}})Q=0$$
 where
\begin{eqnarray}
&&Q=(n^2+3n-2N_2)(1-n\gamma_1^2)
\gamma_1
-((n+2)(n^2\gamma_1^4+1)-6n\gamma_1^2)W^{-1}
\nonumber\\
&&
+2(N_2-1)\gamma_1(1-n\gamma_1^2)W^{-2},
\nonumber\\
&&Q_0=n\bigg(2N_2^4-6(n+1)N_2^3+(n^3+21n+6)N_2^2
-(6n^2+24n+2)N_2+n(n+3)^2\bigg)\gamma_1^2W^2
\nonumber\\
&&
-2n(N_2-1)^2(nN_2-3N_2+n+3)
(1-n\gamma_1^2)\gamma_1W
+(N_2-1)^4(1-n\gamma_1^2)^2
\nonumber\\
&&Q_1=n(nN_2-3N_2+n+3)\gamma_1W
-(N_2-1)^2(1-n\gamma_1^2).
\end{eqnarray} 
If $W=\frac{1-n\gamma_1^2}{(n-1)\gamma_1}$,
then (\ref{equ:beta1}) implies $\beta_1=-\frac{1}{n}$
Then ?? implies $\beta_2=1$.
This is a contradiction. On the other habd we have
$Q_0^2-(N_2-1)^2Q_1^2D_\beta=4\gamma_1^4W^4n^2N_2^2
(N_2-n-1)^6>0$. Hence we must have
\begin{eqnarray}
&&(n^2+3n-2N_2)(1-n\gamma_1^2)
\gamma_1
-((n+2)(n^2\gamma_1^4+1)-6n\gamma_1^2)
W^{-1}
\nonumber\\
&&
+2(N_2-1)\gamma_1(1-n\gamma_1^2)
W^{-2}=0.
\label{equ:WQ}
\end{eqnarray}
(\ref{equ:WP}) and (\ref{equ:WQ})
imply
\begin{eqnarray}
&&(n^2+3n-2N_2)P-2(N_1-1)QW^2
=\bigg((n+2)(n+1)-2(N_1+N_2)\bigg)
\times
\nonumber\\
&&\bigg((n+2-6n\gamma_1^2
+n^2(n+2)\gamma_1^4)W
-(n^2+3n-2)(1-n\gamma_1^2)\gamma_1\bigg)=0.
\end{eqnarray}
If $X$ is not tight, then $N_1+N_2>\frac{(n+2)(n+1)}{2}$
and
we must have
$$(n+2-6n\gamma_1^2
+n^2(n+2)\gamma_1^4)W
-(n^2+3n-2)(1-n\gamma_1^2)\gamma_1=0.$$
This implies
$$W=\frac{(n^2+3n-2)\gamma_1(1-n\gamma_1^2)}
{(2+n-6n\gamma_1^2+2n^2\gamma_1^4+n^3\gamma_1^4)}.
$$ 
Then (\ref{equ:WP})
implies 
$$(1-\gamma_1^2)(1-n^2\gamma_1^2)
(1-n\gamma_1^2)
(n^2\gamma_1^2-n-2
-2n\gamma_1+2n\gamma_1^2)(n^2\gamma_1^2-n-2+2n\gamma_1+2n\gamma_1^2)=0.$$
Hence
$\gamma_1=\frac{1}{n}$ or 
$\gamma_1=\frac{\pm n+\sqrt{n(n+1)(n+4)}}{n(n+2)}
$.
If $\gamma_1=\frac{1}{n}$, then 
$W=1$ and (\ref{equ:alpha1}) implies $\alpha_1=-\frac{1}{n}$
and then $\alpha_2=1$ which is impossible.
On the other hand, if 
$\gamma_1=\frac{\pm n+\sqrt{n(n+1)(n+4)}}{n(n+2)}
$,  then, since $\gamma_1^2$ is a rational number.
Then Proposition \ref{pro:4-3} implies $n=2$. Then we obtain 
 $\gamma_1=\frac{1}{2}$ and 
 again we can introduce a contradiction. 
 Therefore, if $\gamma_1\not=\frac{1}{\sqrt{n}}$, then 
 $(X,w)$ must be a Euclidean tight $4$-design. This completes the proof for Theorem \ref{theo:4-2} (3) and Theorem \ref{theo:1-4}.
 
 \subsection{Proof for Theorem \ref{theo:1-5}}
Assume that a Euclidean $4$-design 
in $\mathbb R^n$ with the parameters
given in Theorem \ref{theo:1-4} (2) (ii) exists. 
Then $X_1\subset S_1=S^{n-1}$ and
$X_2\subset S_2=S^{n-1}(r_2)$.
Let $\boldsymbol y_0=(0,0,\ldots,0,1)\in \mathbb R^{n+1}$.
Let 
$$Y_1=\{(a_1\boldsymbol x,-\frac{1}{2(k-1)})
\mid \boldsymbol x\in X_1\},$$
$$Y_2=\{(a_2\boldsymbol x,\frac{1}{2k})
\mid \boldsymbol x\in X_2\},$$
where $a_1=\frac{\sqrt{(2k-1)(2k-3)}}{2(k-1)}$
and $a_2=\frac{\sqrt{4k^2-1}}{2kr_2}$.
Let $Y=\{\boldsymbol y_0\}\cup Y_1\cup Y_2$.
Then $Y\subset S^n\subset \mathbb R^{n+1}$.
For any $\boldsymbol y_1,\boldsymbol y_2\in Y_1$
$$\boldsymbol y_1\cdot\boldsymbol y_2
=a_1^2\boldsymbol x_1\cdot\boldsymbol x_2+\frac{1}{4(k-1)^2}$$
with $\boldsymbol x_1,\boldsymbol x_2\in X_1$.
Hence 
$\boldsymbol y_1\cdot\boldsymbol y_2=\frac{1}{2k}$
or $-\frac{1}{2(k-1)}$ holds.
If $\boldsymbol y_1,\boldsymbol y_2\in Y_2$, then
$$\boldsymbol y_1\cdot\boldsymbol y_2
=a_2^2\boldsymbol x_1\cdot\boldsymbol x_2
+\frac{1}{4k^2}$$
with $\boldsymbol x_1,\boldsymbol x_2\in X_2$.
Hence $\boldsymbol y_1\cdot\boldsymbol y_2
=\frac{1}{2k}$
or $-\frac{1}{2(k-1)}$ holds. If $\boldsymbol y_1\in Y_1$, $\boldsymbol y_2\in Y_2$
$$\boldsymbol y_1\cdot\boldsymbol y_2
=a_1a_2\boldsymbol x_1\cdot\boldsymbol x_2-\frac{1}{4k(k-1)}$$
with $\boldsymbol x_1\in X_1$ and
$\boldsymbol x_2\in X_2$. Hence $\boldsymbol y_1\cdot\boldsymbol y_2
=\frac{1}{2k}$
or $-\frac{1}{2(k-1)}$ holds.
Thus we obtain a $2$-distance set $Y$
on $S^n$ whose cardinality attains the
Fisher bound, $\frac{(n+2)(n+1)}{2}+1$, for the $2$-distance set 
on the unit sphere $S^n$.
Therefore $Y$ is a spherical tight $4$-design on
$S^n$.
Conversely,
assume a spherical tight $4$-design $Y$ exists
on $S^n$. Then there exists an integer $k\geq 2$ satisfying
$n+4=(2k-1)^2$ (see \cite{B-D-1,B-D-2}).
It is known that if $Y$ is a spherical $4$-design
then
$Y^\tau$ is a spherical $4$-design,
where $Y^\tau$ is the image of $Y$ under 
$\tau\in O(n+1)$ (orthogonal group of degree $n+1$). Hence we may assume the unit vector 
$\boldsymbol y_0=(0,0,\ldots,0,1)$
is contained in $Y$. 
It is also known that $Y$ is a $2$-distance set
and $\boldsymbol y_1\cdot \boldsymbol y_2
=-\frac{1}{2(k-1)}$ or $\boldsymbol y_1\cdot \boldsymbol y_2
=\frac{1}{2k}$ holds for any $\boldsymbol{\boldsymbol y_1,\boldsymbol y_2}\in Y$ (see \cite{D-G-S}, 
also Appendix I).
Let 
$$Y_1=\left\{\boldsymbol y\in Y\ \bigg|\ \boldsymbol y\cdot\boldsymbol y_0
=-\frac{1}{2(k-1)}\right\},\quad
Y_2=\left\{\boldsymbol y\in Y\ \bigg|\ \boldsymbol y\cdot\boldsymbol y_0
=\frac{1}{2k}\right\}.$$
It is also known that $|Y_1|=2(2k+1)(k-1)^3$ and
$|Y_2|=2(2k-3)k^3$ holds.
$(n+1)$st coordinate of the vectors in $Y_1$ are
$-\frac{1}{2(k-1)}$ and  $Y_2$ are $\frac{1}{2k}$
respectively.
Let 
$$X_1=\{\boldsymbol x\in \mathbb R^n
\mid (\boldsymbol x,-\frac{1}{2(k-1)})\in Y_1\},
\quad X_2=\{\boldsymbol x\in \mathbb R^n
\mid (\boldsymbol x,\frac{1}{2k})\in Y_2\}.$$
$X_1$ is on the sphere of radius $r_1=\frac{\sqrt{(2k-1)(2k-3)}}
{2(k-1)}$ and
$X_2$ is on the sphere of radius
$r_2=\frac{\sqrt{(2k-1)(2k+1)}}{2k}$.
Then we obtain
\begin{eqnarray}&&A(X_1,X_1)=\left\{\frac{k-2}{k(2k-3)}
,-\frac{1}{2k-3}
\right\},\nonumber\\
&&A(X_2,X_2)=\left\{\frac{1}{2k+1},
-\frac{k+1}{(k-1)(2k+1)}\right\},
\nonumber\\
&&A(X_1,X_2)=\left\{\frac{1}{\sqrt{(2k+1)(2k-3)}}
\left(=\frac{1}{\sqrt{n}}\right),
\frac{-1}{\sqrt{(2k+1)(2k-3)}}
\left(=\frac{-1}{\sqrt{n}}\right)\right\}.
\nonumber
\end{eqnarray}
Since $Y$ is a spherical tight 4-design on $S^n$,
$Y$ has the structure of a Q-polynomial association scheme.
Also it is proved that the Q-polynomial association scheme is three regular. Hence 
$X=X_1\cup X_2$ has the structure of a coherent configuration with the same parameter as given in
Theorem \ref{theo:1-4} (2) (ii).

\section{Theorem \ref{theo:1-6}}
In this section, we consider the case 
when $(X,w)$ is Euclidean tight $4$-design
and $N_1\geq n+2$. 
Hence $N_2=\frac{(n+2)(n+1)}{2}-N_1\geq n+2$.
Then (\ref{equ:WP}) and (\ref{equ:WQ})
imply
\begin{eqnarray}
&&W=\frac{(n+2)(n^2\gamma_1^4+1)-6n\gamma_1^2
+\varepsilon\sqrt{D_1}
}{4(N_1-1)\gamma_1(1-n\gamma_1^2)}
\label{equ:W1}\\
&&W^{-1}=\frac{(n+2)(n^2\gamma_1^4+1)-6n\gamma_1^2
-\varepsilon\sqrt{D_2}
}{4(N_2-1)\gamma_1(1-n\gamma_1^2)}
\nonumber
\end{eqnarray}
where $\varepsilon=\pm1$ and
\begin{eqnarray}&&D_i
=16N_i^2\gamma_1^2(1-n\gamma_1^2)^2
-8(n+2)(n+1)N_i\gamma_1^2(1-n\gamma_1^2)^2
+(n+2)^2\nonumber\\
&&-4n^2\gamma_1^2+2n^2(n^2-4n-2)\gamma_1^4
-4n^4\gamma_1^6+n^4(n+2)^2\gamma_1^8.
\end{eqnarray}
By (\ref{equ:gamma-m}) we have
$\gamma_1=\sqrt{\frac{N_1-m_1}{nm_1}}
=\sqrt{\frac{N_2-m_2}{nm_2}}$
and $\frac{N_1}{N_2}=\frac{m_1}{m_2}$.
Since $0<\gamma_1\neq \frac{1}{\sqrt{n}}\leq 1$,
we also have $\frac{N_i}{2}< m_i \leq N_i-1$, $i=1,2$.
Since 
$D_i\geq 0$, we have
$$\frac{N_i}{2}< m_i\leq \frac{N_i}{2}+ \frac{N_i}{2}\sqrt{\frac{K_2(n,N_i)}{K_1(n,N_i)}},$$
$$\frac{N_i}{2}+ \frac{N_i}{2}\sqrt{\frac{K_3(n,N_i)}{K_1(n,N_i)}}\leq m_i<N_i$$
where
$$K_1(n,x)=8x(n^2+3n+2-2x)+n^3+2n^2+n,$$
\begin{eqnarray}K_2(n,x)=
4x(n^2+3n+2-2x)-n(n+1)(n+7)\nonumber\\
-4\sqrt{(N_i-1)
(3n+n^2-2x)(x-n-1)(n+n^2-2x)}
\end{eqnarray}
and
\begin{eqnarray}K_3(n,x)=
4x(n^2+3n+2-2x)-n(n+1)(n+7)\nonumber\\
+4\sqrt{(x-1)
(3n+n^2-2x)(x-1-n)(n+n^2-2x)}
\end{eqnarray}

Then we can express $\alpha_1$ and $\alpha_2$
interms of $n,\ N_1$ and $m_1$,  and
$\beta_1$ and $\beta_2$ 
interms of $n, N_1$ and $m_1$.
Then we can express 
the ratio $\left(\frac{2-\alpha_1-\alpha_2}{\alpha_1-\alpha_2}\right)^2$
interms of $n,N_1$ and $m_1$ and 
$\left(\frac{2-\beta_1-\beta_2}{\beta_1-\beta_2}\right)^2$
interms of $n,N_2$ and $m_2$.
Then we have
\begin{eqnarray}
&&\left(\frac{2-\alpha_1-\alpha_2}{\alpha_1-\alpha_2}\right)^2=
F_\varepsilon(n,N_1,m_1)\\
&&\left(\frac{2-\beta_1-\beta_2}{\beta_1-\beta_2}\right)^2=
F_\varepsilon(n,N_2,m_2),
\end{eqnarray}
where
\begin{equation}
F_{\varepsilon}(n,x,y)=\frac{F_1(n,x,y)-4\varepsilon F_2(n,x,y)\sqrt{n
F_4(n,x,y)}}{F_3(n,x,y)}
\end{equation}
and
\begin{eqnarray}
&&F_1(n,x,y)=8n(n+1)^2(n+1+4x)y^4
-16xn(n+1)^2(n+1+4x)y^3
\nonumber\\
&&\qquad+4x^2(8n^3+24xn+4x^2-8n-4x^2n-4x+3n^4+30n^2x-3n^2+10n^3x)y^2
\nonumber\\
&&\qquad-4x^3(2n^3-10n-4x^2n+16xn-9n^2+4x^2+2n^3x-4x+n^4+14n^2x)y
\nonumber\\
&&\qquad+x^3(-1+x)(4x^2+4xn+6n^3+n^4+9n^2)\nonumber\\
&&F_2(n,x,y)=((n+1)^2y^2-x(n+1)^2y+x^2(-1+x))\nonumber\\
&&F_3(n,x,y)=16n(n+1)^3y^4-32xn(n+1)^3y^3
\nonumber\\
&&
\qquad+8x^2n(-2+6x+5n-2xn+10n^2+3n^3)y^2
\nonumber\\
&&
\qquad-8x^3n(6x-4-2xn-n+4n^2+n^3)y
\nonumber\\
&&
\qquad+x^3(-1+x)(4x^2-4xn-4n^2x+6n^3+n^4+9n^2)
\nonumber\\
&&F_4(n,x,y)=(4n+8n^2+4n^3+64x-64x^2+32n^2x+96xn)y^4
\nonumber\\
&&
\qquad-8x(n+2n^2+n^3+16x-16x^2+8n^2x+24xn)y^3
\nonumber\\
&&
\qquad+4x^2(20x-20x^2+5n+30xn+7n^2+10n^2x+2n^3)y^2
\nonumber\\
&&
\qquad-4x^3(-4x^2+4x+6xn+4n+2n^2x+5n^2+n^3)y
\nonumber\\
&&
\qquad+nx^4(n+2)^2.
\end{eqnarray} 

Since numerical experiments for
small $n$ shows that every integral conditions
are satisfied only if $n+3$ is the
square of an odd integer.
In that case (if integral conditions are satisfied), 
numerical experiments for
small $n$ shows
$\left(\frac{2-\alpha_1-\alpha_2}{\alpha_1-\alpha_2}
\right)^2=\left(\frac{2-\beta_1-\beta_2}{\beta_1-\beta_2}
\right)^2=n+3$.
In the following we
assume that $n+3$ is a square of an odd integer.

$\left(\frac{2-\beta_1-\beta_2}{\beta_1-\beta_2}
\right)^2=n+3$. Then, for both $\varepsilon=1$ and
$\varepsilon=-1$, we obtain

\begin{equation}m_2=\frac{N_2}{2}+
\frac{1}{2}\sqrt{\frac{-(n-3)N_1N_2^2
+(N_2-N_1)
\sqrt{2n(n+3)N_1N_2^3}}
{(n+1)^2N_1}}
\label{equ:m2}
\end{equation}
and
\begin{equation}m_1=\frac{N_1}{2}+\frac{1}{2}
\sqrt{\frac{-(n-3)N_1^2N_2
+(N_2-N_1)\sqrt{2n(n+3)N_1^3N_2
})}{(n+1)^2N_2}}
\label{equ:m1}
\end{equation}
(Note that $N_2=\frac{(n+2)(n+1)}{2}-N_1$.)
Conversely,
if we assume $\left(\frac{2-\alpha_1-\alpha_2}{\alpha_1-\alpha_2}
\right)^2=n+3$, then we obtain 
(\ref{equ:m1}), 
(\ref{equ:m2}) for both $\varepsilon =1$
and $\varepsilon =-1$.

Since $m_1$ and $m_2$ are integers,
satisfying $m_i>\frac{N_i}{2}$
we must have
$N_1< \frac{(n+2)(n+3)}{6}$
(equivalently $N_2> \frac{n(n+2)}{3}
$).
We note that if $N_1= \frac{(n+2)(n+3)}{6}$,
then $m_1=\frac{N_1}{2}$ and contradicts
the fact $\gamma_1<\frac{1}{\sqrt{n}}$.

Numerical experiments suggest us that $N_1=\frac{n(n+1)}{6}$
($N_2=\frac{(n+3)(n+1)}{3}$) 
gives good conditions.
Actually, for any $n$ satisfying 
$n+3=(6k-3)^2$, with an integer $k\geq 2$,
let  

$$N_1 =\frac{n(n+1)}{6}= (6k^2-6k+1)(36k^2-36k+7)
$$
and
$$m_1 =\frac{n(n+4)}{12}= (6k^2-6k+1)(18k^2-18k+5).$$
Then we obtain
\begin{eqnarray}
&&N_2 =\frac{(n+3)(n+1)}{3}= 3(36k^2-36k+7)(2k-1)^2,\nonumber\\
&&
m_2=\frac{(n+3)(n+4)}{6} = 3(2k-1)^2(18k^2-18k+5).
\nonumber\\
&&\gamma_1=\sqrt{\frac{9k^2-9k+1}
{(18k^2-18k+3)(18k^2-18k+5)}}
\end{eqnarray}
Moreover
let $\varepsilon =-1$ in equation (\ref{equ:W1}). Then we obtain
$$W=\frac{\sqrt{(18k^2-18k+5)^3(18k^2-18k+3)}}
{9(2k-1)^2\sqrt{(9k^2-9k+1)^3}
}
$$
and
$\left(\frac{2-\alpha_1-\alpha_2}{\alpha_1-\alpha_2}\right)^2
=\left(\frac{2-\beta_1-\beta_2}{\beta_1-\beta_2}\right)^2=n+3=(6k-3)^2$.
We also have
$$A(X_1,X_1)=
\left\{\frac{18k^2-27k+8}{6(9k^2-9k+1)(2k-1)},\
\frac{18k^2-9k-1}{6(9k^2-9k+1)(2k-1)}
\right\}
$$
$$A(X_2,X_2)
=\left\{\frac{36k^3-54k^2+25k-4}{2(6k^2-6k+1)(18k^2-18k+5)},\
\frac{36k^3-54k^2+25k-3}{2(6k^2-6k+1)(18k^2-18k+5)}
\right\}.
$$
$$r=\frac{\sqrt{(18k^2-18k+5)(18k^2-18k+3)}}
{\sqrt{9k^2-9k+1}}.
$$
Thus we can determined all the parameters in terms of $k$ and also we can express all the possible intersection numbers
of the corresponding coherent configuration
in polynomials of $k$. The reader can find them 
in Appendix II. 

Exhaustive numerical experiments 
for the case $\varepsilon=-1$ in (\ref{equ:W1}) for
every $n$ up to $n=222$, shows that there is no feasible parameter other than
this family.
Also
exhaustive numerical experiments 
for the case $\varepsilon=1$ in (\ref{equ:W1}) for
every $n$ up to $n=222$, shows that
only $n=22=(6-1)^2-3$, $N_1=33$, 
$N_2=243$, $m_1=22$, $m_2=162$ satisfies 
every requirement for Euclidean tight 4-design
on 2 concentric spheres.
The Euclidean tight 4-design with this parameter is
constructed and unique (Theorem III in \cite{B-2}).
As for this parameter we can consider from 
a different view point explained in the next section.

\section{An additional remark}
We also proved the 
following.
\begin{theo}\label{theo:4-4} If $\alpha_1=0$, then
$(X,w)$ is similar to one of the following Euclidean tight 
$4$-designs.

(1) $n=4$, $N_1=6$, $N_2=9$ and given in Theorem II in \cite{B-2}.

(2) $n=22$, $N_1=33$, $N_2=243$ and given in Theorem III in \cite{B-2}..
\end{theo}
{\bf Proof}
Let $\alpha_1=0$. Then (\ref{equ:alpha2})
implies $\alpha_2=-\frac{N_1-n}{n}$.
Then we obtain
$$\frac{2-\alpha_1-\alpha_2}{\alpha_1-\alpha_2}
=\frac{N_1+n}{N_1-n}=1+\frac{2n}{N_1-n}.$$
Since $\frac{2n}{N_1-n}$ must be an even integer,
$\frac{n}{N_1-n}$ is a positive integer. Hence
$N_1\leq 2n$ and $n$ is a multiple of
$N_1-n$.
Since 
$p_{\alpha_1,\alpha_2}^{\alpha_2}=
p_{\alpha_2,\alpha_1}^{\alpha_2}=
\frac{N_1^2\alpha_1(n\alpha_1+1)
(N_1-n-1)}{\left(2n\alpha_1-n+N_1+
n(N_1-1)\alpha_1^2\right)^2}=0$ (see Appendix III of this paper),
$X_1$ is a union of $1$-distance set containing $p_{\alpha_2,\alpha_2}^{\alpha_0}
+1=\frac{n}{N_1-n}+1$ points and mutually perpendicular
to each other.
Let $\frac{n}{N_1-n}+1=q$.
Then $X_1=\Delta_1\cup\cdots
\cup\Delta_d$, $|\Delta_i|=q$
and $N_1=dq$. Then $n=d(q-1)$
and each $\Delta_i$ is a $1$-distance set
in $\mathbb R^{q-1}$. 
Since $\alpha_2=-\frac{1}{q-1}$,
$\Delta_i$, $(1\leq i\leq d)$ is a regular simplex in 
$\mathbb R^{q-1}$. 
On the other hand (\ref{equ:V}) and 
(\ref{equ:II}) imply
\begin{eqnarray}
&&\gamma_1^2=
\frac{1}{n^3}\left(-n^2+N_1(3n-N_1)
\pm \sqrt{N_1(N_1-n)(2n-N_1)(3n-N_1)}\right)
\nonumber\\
&&
=\frac{q^2-q-1\pm\sqrt{q(q-2)(2q-3)}}{d(q-1)^3}.
\end{eqnarray}
Since $\gamma_1^2$ is a rational number
$q(q-2)(2q-3)$ must be the square of an integer.

\begin{pro}\label{pro:Kaneko2}
$q(q-2)(2q-3)$ is the square of an integer
if and only if $q=2$ and $q=3$.
\end{pro}
{\bf Proof} Kaneko \cite{K}.\\

If $q=2$, then $d=n$ and we obtain
$\gamma_1=\frac{1}{\sqrt{n}}$. This contradicts the
assumption $\gamma_1\not=\frac{1}{\sqrt{n}}$.
If $q=3$, then $N_1=3d$, $n=2d$,
$N_2=2d^2+1$ and
$\gamma_1=\frac{1}{\sqrt{d}}$ or 
$\gamma_1=\frac{1}{2\sqrt{d}}$.
Then
$\gamma_1=\frac{1}{\sqrt{d}}$ and 
(\ref{equ:I})
implies $r_2^3(2d^2+1)+3d\sqrt{d}=0$.
This is a contradiction.
Hence $\gamma_1=\frac{1}{2\sqrt{d}}$
holds.
Then (\ref{equ:I}) and (\ref{equ:II})
imply $r_2=\sqrt{d}$, $w_2=\frac{3}{2d^2+1}$.
Then (\ref{equ:III})
and (3-i-c) implies
$$\beta_1=\frac{-1+\sqrt{8d-7}}{4d},\quad
\beta_2=\frac{-1-\sqrt{8d-7}}{4d}.$$
Since $\beta_1$ and $\beta_2$ are rational numbers
$8d-7=(2k-1)^2$ with a integer $k\geq 1$.
Then
$$\frac{2-\beta_1-\beta_2}{\beta_1-\beta_2}
=k-\frac{k-5}{2k-1}.$$
If $k\geq 6$, then $0<\frac{k-5}{2k-1}<1$
and $\frac{2-\beta_1-\beta_2}{\beta_1-\beta_2}$
cannot be an integer.
Hence we must have 
$k=1,2$ or $5$. Since $N_1\geq n+2$, we have
$k=2$, $d=2$, $n=4$, $N_1=6$, $N_2=9$;
and $k=5$, $d=11$ $n=22$, $N_1=33$,
$N_1=243$.
Theorem II and  Theorem III in \cite{B-2} implies 
Theorem \ref{theo:4-4}. \\

\section{Appendix I}
Feasible parameters of the Euclidean $4$-design
$(X,w)$ given in Theorem \ref{theo:1-4}
(2) (ii) and the intersection numbers of the coresponding
coherent configuration.\\

\noindent
\qquad$n=(2k-1)^2-4$,\\ 
\qquad$|X_1|=2(2k+1)(k-1)^3$, 
\qquad$|X_2|=2k^3(2k-3)$,\\
\qquad$A(X_1,X_1)=\{\frac{k-2}{k(2k-3)},\
-\frac{1}{2k-3}\}$, 
\qquad$A(X_2,X_2)=\{\frac{1}{2k+1},
-\frac{k+1}{(k-1)(2k+1)}\}$,\\
\qquad$A(X_1,X_2)=\{\frac{1}{\sqrt{n}},\
-\frac{1}{\sqrt{n}}\}$,\\ 
\qquad$r_1=1$, \qquad$w_1=1$,
\qquad$w_2= \frac{(2k+1)^2(k-1)^4}{(2k-3)^2k^4}r_2^{-4}$,\\

\noindent
{\bf Intersection matrices and Character tables of the association scheme for
$X_1$}\\

$B_1^{(1)}=\left[
\begin{array}{ccc}0& 1& 0\\ 
k^3(2k-3)& (k+1)(k^2-k-1)k& (k-1)k^3\\
 0& (k^2-k-1)(k-1)^2& k^3(k-2)
 \end{array}\right]$,\\
 
 $B_2^{(1)}=\left[\begin{array}{cc}
 0&0\\
 0& (k^2-k-1)(k-1)^2\\
 (k-1)(2k-3)(k^2-k-1)& (k-2)(k-1)(k^2-k-1)
 \end{array}\right.$,
 
\hfill $\left.\begin{array}{c}
 1\\
 k^3(k-2)\\
 (k-1)(k-2)(k^2-2k-2)
 \end{array}\right]$,\\

$P_1 =\left[
\begin{array}{ccc} 1&k^3(2k-3)&(k-1)(2k-3)(k^2-k-1)
\\
1& k^2(k-2)&-1-k^2(k-2)\\
1&-k&-1+k \end{array}\right]$,

$Q_1=\left[
\begin{array}{ccc} 1& (2k+1)(2k-3)& 2(2k-3)(k^2-k-1)k\\
1&\frac{ (k-2)(2k+1)}{k}& -\frac{2(k^2-k-1)}{k}\\
1& -2k-1& 2k\end{array}\right]$,\\

\noindent
{\bf Intersection matrices and Character tables of the association scheme for
$X_2$}\\

$B_1^{(2)}=\left[
\begin{array}{ccc} 0& 1& 0\\
(2k+1)(k^2-k-1)k& (k+1)(k^2-3)k& (k+1)(k^2-k-1)k\\
0& (k+1)(k-1)^3& (k^2-k-1)k^2
\end{array}\right]$,

$B_2^{(2)}=\left[
\begin{array}{ccc} 0& 0& 1\\
0& (k+1)(k-1)^3& (k^2-k-1)k^2\\
(2k+1)(k-1)^3& (k-1)^3k& (k-2)(k-1)(k^2-k-1)
\end{array}\right]$

$P_2=\left[
\begin{array}{ccc} 1& (2k+1)(k^2-k-1)k& (k-1)(2k^3-3k^2+1)\\
1& k(k^2-k-1)& -(k-1)(k^2-1)\\
1& -k& k-1]])
\end{array}\right]$

$Q_2 =\left[
\begin{array}{ccc} 
1& (2k+1)(2k-3)& 2(k-1)(2k+1)(k^2-k-1)\\
1& 2k-3& -2k+2\\
1& -\frac{(2k-3)(k+1)}{k-1}& \frac{2(k^2-k-1)}{k-1}\end{array}\right]$\\

$p_{\gamma_1,\gamma_1}^{\alpha_0}=k^3(2k-3)
$,\quad
$p_{\gamma_1,\gamma_1}^{\beta_0}= (2k+1)(k-1)^3
$
\vskip.2cm
$p_{\gamma_2,\gamma_2}^{\alpha_1}=(k^2-k-1)k^2
$,\quad
$p_{\gamma_1,\gamma_2}^{\alpha_1}=(k-1)^2k^2
$,\quad
$p_{\gamma_1,\gamma_1}^{\alpha_1}=(k^2-k-1)k^2
$
\vskip.2cm
$p_{\gamma_2,\gamma_2}^{\alpha_2}=k^3(k-2)
$,\quad
$p_{\gamma_1,\gamma_2}^{\alpha_2}=(k-1)k^3
$,\quad
$p_{\gamma_1,\gamma_1}^{\alpha_2}=k^3(k-2)
$
\vskip.2cm
$p_{\gamma_1,\gamma_2}^{\beta_1}= (k-1)^3k
$,\quad
$ p_{\gamma_2,\gamma_2}^{\beta_1}= (k+1)(k-1)^3$,\quad
$ p_{\gamma_1,\gamma_1}^{\beta_1}= (k+1)(k-1)^3
$
\vskip.2cm
$p_{\gamma_1,\gamma_2}^{\beta_2}=(k-1)^2k^2
$,\quad
$p_{\gamma_2,\gamma_2}^{\beta_2}= (k^2-k-1)(k-1)^2
$,\quad
$ p_{\gamma_1,\gamma_1}^{\beta_2}=(k^2-k-1)(k-1)^2
$,\\
\vskip.2cm
$ p_{\gamma_2,\beta_2}^{\gamma_1}= (k-1)^2k^2
$,\quad
$ p_{\alpha_2,\gamma_2}^{\gamma_1}=(k^2-k-1)(k-1)^2
$,\quad
$ p_{\gamma_1,\beta_1}^{\gamma_1}=(k+1)(k^2-k-1)k
$,\\

$ p_{\gamma_2,\beta_1}^{\gamma_1}= (k^2-k-1)k^2
$
,\quad
$ p_{\alpha_1,\gamma_1}^{\gamma_1}=(k^2-k-1)k^2
$,\quad
$ p_{\alpha_1,\gamma_2}^{\gamma_1}=(k-1)^2k^2
$,\\

$ p_{\gamma_1,\beta_2}^{\gamma_1}= (k^2-k-1)(k-1)^2
$,\quad
$ p_{\alpha_2,\gamma_1}^{\gamma_1}=(k-2)(k-1)(k^2-k-1)
 $,\\
\vskip.2cm
$ p_{\gamma_2,\beta_2}^{\gamma_2}=(k^2-k-1)(k-1)^2
$
,\quad
$ p_{\alpha_2,\gamma_2}^{\gamma_2}=(k-2)(k-1)(k^2-k-1)
$,\\

$ p_{\gamma_1,\beta_2}^{\gamma_2}=(k-1)^2k^2
$,
\quad
$ p_{\gamma_1,\beta_1}^{\gamma_2}=(k^2-k-1)k^2
$,
\quad
$ p_{\alpha_1,\gamma_2}^{\gamma_2}=(k^2-k-1)k^2
$,\\

$ p_{\alpha_2,\gamma_1}^{\gamma_2}=(k^2-k-1)(k-1)^2
$,
\quad
$ p_{\gamma_2,\beta_1}^{\gamma_2}=(k+1)(k^2-k-1)k
$,\quad
$ p_{\alpha_1,\gamma_1}^{\gamma_2}=(k-1)^2k^2
$.\\

\noindent
In above $p_{a,b}^{c}
=p_{b,a}^{c}$ holds for any $a,b,c\in
\{\alpha_i,\beta_j, \gamma_k\mid
i,\ j=0,1,2, k=1,2\}$.\\

\noindent
{\bf Sperical tight 4-design on $S^n\subset \mathbb R^{n+1}$}

If sperical tight 4-design $Y\subset S^n\subset\mathbb R^{n+1}$
exists, then we must have $n+4=(2k-1)^2$
with an integer $k\geq 2$.
Then $A(Y)=\{
\frac{-1-\sqrt{n+4}}{n+3},\frac{-1+\sqrt{n+4}}{n+3}
\}=\{-\frac{1}{2(k-1)}, \frac{1}{2k}\}$.
Then $Y$ has the structure of an association scheme
whose second eigen-matrix is given by
$$\left[\begin{array}{lll}
Q_{n+1,0}(1)&Q_{n+1,1}(1)&Q_{n+1,2}(1)\\
Q_{n+1,0}(-\frac{1}{2(k-1)})
&Q_{n+1,1}(-\frac{1}{2(k-1)})
&Q_{n+1,2}(-\frac{1}{2(k-1)})
\\
Q_{n+1,0}(\frac{1}{2k})&
Q_{n+1,1}(\frac{1}{2k})&
Q_{n+1,2}(\frac{1}{2k})
\end{array}
\right]$$
$$=\left[\begin{array}{ccc}
1& 4k^2-4k-2& 2(2k+1)(2k-3)k(k-1)\\
1& -\frac{2k^2-2k-1}{k-1}& \frac{(2k-3)k}{k-1}
\\
1& \frac{2k^2-2k-1}{k}& -\frac{(2k+1)(k-1)}{k}
\end{array}
\right]
.$$
This indicate that $E_1$ induces the projection 
of the association scheme into the unit sphere
$S^n\subset\mathbb R^{n+1}$. 
The character table of $Y$
is given by
$$\left[
\begin{array}{ccc}
1&2(2k+1)(k-1)^3&2(2k-3)k^3\\ 
1&-(2k+1)(k-1)^2&(2k-3)k^2\\
1&k-1&-k
\end{array}
\right].
$$
Let $\boldsymbol u_0$ be a fixed point in $Y$.
We may assume 
$ \boldsymbol u_0=(0,0,\ldots,0,1)$.
Let $Y_2=\{\boldsymbol y\in Y\mid \boldsymbol u_0\cdot
\boldsymbol y=\frac{1}{2k}\}$ and
$Y_1=\{\boldsymbol y\in Y\mid \boldsymbol u_0\cdot
\boldsymbol y=-\frac{1}{2(k-1)}\}$.
$|Y_2|=2(2k-3)k^3$ and $|Y_1|=2(2k+1)(k-1)^3$.

\section{Appendix II}
The feasible parameters of the Euclidean tight $4$-design given in Theorem \ref{theo:1-6} and
intersection numbers of the corresponding 
coherent configuration.\\

$n=(6k-3)^2-3$, 

$|X_1|=(6k^2-6k+1)(36k^2-36k+7)$,\quad
$|X_2|=3(36k^2-36k+7)(2k-1)^2$,

$A(X_1,X_1)=${\Large$\left\{
\frac{18k^2-27k+8}{6(9k^2-9k+1)(2k-1)},\ 
-\frac{18k^2-9k-1}{6(9k^2-9k+1)(2k-1)}
\right\}$,}\\

$A(X_2,X_2)=$\\
\qquad{\Large$\left\{\frac{36k^3-54k^2+25k-4}
{2(6k^2-6k+1)(18k^2-18k+5)},\
-\frac{36k^3-54k^2+25k-3}
{2(6k^2-6k+1)(18k^2-18k+5)}\right\}$},\\
 
 $A(X_1,X_2)=$\\
\qquad {\Large$\left\{\begin{array}{l}\sqrt{\frac{36k^2-36k+4}{(36k^2-36k+6)
 (36k^2-36k+10)}},
 -\sqrt{\frac{36k^2-36k+10}{(36k^2-36k+6)
 (36k^2-36k+4)}}
 \end{array}\right\}$},\\

 $r_1=1$,
 \quad $r_2= \sqrt{\frac{3(18k^2-18k+5)(6k^2-6k+1)}
 {9k^2-9k+1}}
$, \\

$w_1=1$, \quad $w_2= \frac{1}{81(2k-1)^4}$.\\

\noindent
{\bf Intersection matrices and the Character tables of the association scheme for
$X_1$}\\

$B_1^{(1)}=\left[
\begin{array}{cc}0& 1\\
6(-1+2k)(9k^2-9k+1)k& 54k^4-45k^3-12k^2+7k+1\\
0& (3k-2)(k-1)(18k^2-9k-1)
 \end{array}\right.$,\\
 
 \hfill$\left.\begin{array}{c}
 0\\
 (18k^2-9k-1)k(3k-2)\\
 k(3k-1)(18k^2-27k+8)
 \end{array}\right]$\\

$B_1^{(2)}=\left[
\begin{array}{cc}0& 0\\
0& (3k-2)(k-1)(18k^2-9k-1)\\
6(k-1)(-1+2k)(9k^2-9k+1)& (18k^2-27k+8)(k-1)(3k-1)
 \end{array}\right.$\\
 
 \hfill
$\left. \begin{array}{c} 
1\\
 k(3k-1)(18k^2-27k+8)\\
 54k^4-171k^3+177k^2-64k+5
\end{array}\right]$,\\

\vskip.5cm
$P_1=\left[
\begin{array}{ccc} 
1& 6(-1+2k)(9k^2-9k+1)k& 6(k-1)(-1+2k)(9k^2-9k+1)\\1& -3k+1& 3k-2\\1& k(18k^2-27k+8)& -(k-1)(18k^2-9k-1)
 \end{array}\right]$,\\

$Q_1=\left[
\begin{array}{ccc} 1& 6(36k^2-36k+7)(k-1)k& 36k^2-36k+6\\
1& -\frac{(3k-1)(k-1)(36k^2-36k+7)}{(-1+2k)(9k^2-9k+1)}& 
\frac{(18k^2-27k+8)(6k^2-6k+1)}{(-1+2k)(9k^2-9k+1)}\\1&\frac{ k(3k-2)(36k^2-36k+7)}{(-1+2k)(9k^2-9k+1)}& -\frac{(18k^2-9k-1)(6k^2-6k+1)}
{(-1+2k)(9k^2-9k+1)} \end{array}\right]$,\\

\noindent
{\bf Intersection matrices and the Character tables of the association scheme for
$X_2$}\\

$B_2^{(1)}=\left[
\begin{array}{cc} 0& 1\\
2(6k^2-6k+1)(18k^2-18k+5)& (9k^2-9k+1)(12k^2-10k+3)\\
0& (3k-2)(36k^3-54k^2+25k-3)\end{array}\right.$,\\

\hfill
$\left. \begin{array}{c} 
0\\
 (3k-2)(36k^3-54k^2+25k-3)\\
 (36k^3-54k^2+25k-4)(3k-1) 
\end{array}\right]$,\\

$B_2^{(2)}=\left[
\begin{array}{cc} 0& 1\\
0& (3k-2)(36k^3-54k^2+25k-3)\\
2(6k^2-6k+1)(18k^2-18k+5)& (36k^3-54k^2+25k-4)(3k-1)\end{array}\right.$,\\

\hfill
$\left. \begin{array}{c} 0\\
(36k^3-54k^2+25k-4)(3k-1)\\
 (9k^2-9k+1)(12k^2-14k+5)
\end{array}\right]$,\\

\vskip.2cm

$P_2=\left[
\begin{array}{ccc} 
1& 2(6k^2-6k+1)(18k^2-18k+5)& 2(6k^2-6k+1)(18k^2-18k+5)\\
1& -3k+1& 3k-2\\
1& 36k^3-54k^2+25k-4& 3-36k^3+54k^2-25k
\end{array}\right]$,\\

\vskip.2cm
$Q_2 =\left[
\begin{array}{ccc} 
1& 2(6k^2-6k+1)(36k^2-36k+7)& 36k^2-36k+6\\
1& -\frac{(3k-1)(36k^2-36k+7)}{18k^2-18k+5}&
\frac{ 3(36k^3-54k^2+25k-4)}{18k^2-18k+5}\\
1& \frac{(3k-2)(36k^2-36k+7)}{18k^2-18k+5}& 
-\frac{3(36k^3-54k^2+25k-3)}{18k^2-18k+5} \end{array}\right]$,\\
\vskip.2cm
$p_{\gamma_1,\gamma_1}^{\alpha_0}=3(18k^2-18k+5)(2k-1)^2
$,\quad
$p_{\gamma_1,\gamma_1}^{\beta_0}=(6k^2-6k+1)(18k^2-18k+5)
$,

\vskip.2cm

$p_{\gamma_2,\gamma_2}^{\alpha_1}=(2k-1)(54k^3-72k^2+15k+4)
$,
\quad
$p_{\gamma_1,\gamma_2}^{\alpha_1}=(3k-2)(2k-1)(18k^2-18k+5)
$,

$p_{\gamma_1,\gamma_1}^{\alpha_1}=(2k-1)(3k-1)(18k^2-18k+5)
$,
\vskip.2cm
$p_{\gamma_2,\gamma_2}^{\alpha_2}=(54k^3-90k^2+33k-1)(2k-1)
$,
\quad
$p_{\gamma_1,\gamma_2}^{\alpha_2}= (2k-1)(3k-1)(18k^2-18k+5)
$,

$p_{\gamma_1,\gamma_1}^{\alpha_2}=(3k-2)(2k-1)(18k^2-18k+5)
$,
\vskip.2cm
$p_{\gamma_1,\gamma_2}^{\beta_1}= (2k-1)(3k-2)(9k^2-9k+1)
$,
\quad
$ p_{\gamma_2,\gamma_2}^{\beta_1}= (9k^2-9k+1)k(6k-5)
$,

$ p_{\gamma_1,\gamma_1}^{\beta_1}= (3k-1)(18k^3-27k^2+14k-3)$,
\vskip.2cm

$p_{\gamma_1,\gamma_2}^{\beta_2}= (3k-1)(9k^2-9k+1)(2k-1)
$,
\quad
$p_{\gamma_2,\gamma_2}^{\beta_2}= (9k^2-9k+1)(6k-1)(k-1)
$,

$ p_{\gamma_1,\gamma_1}^{\beta_2}= (3k-2)(18k^3-27k^2+14k-2)
$,
\vskip.2cm

$ p_{\gamma_2,\beta_2}^{\gamma_1}=2(3k-1)(9k^2-9k+1)(2k-1)
$,
\quad
$ p_{\alpha_2,\gamma_2}^{\gamma_1}=2(3k-1)(k-1)(9k^2-9k+1)
$,

$ p_{\gamma_1,\beta_1}^{\gamma_1}=2(3k-1)(18k^3-27k^2+14k-3)
$,
\quad
$ p_{\gamma_2,\beta_1}^{\gamma_1}=2(2k-1)(3k-2)(9k^2-9k+1)
$,

$ p_{\alpha_1,\gamma_1}^{\gamma_1}=2k(3k-1)(9k^2-9k+1)
$,
\qquad\qquad
$ p_{\alpha_1,\gamma_2}^{\gamma_1}=2k(3k-2)(9k^2-9k+1)
$,

$ p_{\gamma_1,\beta_2}^{\gamma_1}=2(3k-2)(18k^3-27k^2+14k-2)
$,
\quad
$ p_{\alpha_2,\gamma_1}^{\gamma_1}= 2(k-1)(9k^2-9k+1)(3k-2)
$,

\vskip.2cm
$ p_{\gamma_2,\beta_2}^{\gamma_2}=(6k-1)(k-1)(18k^2-18k+5)
$,
\quad
$ p_{\alpha_2,\gamma_2}^{\gamma_2}=(k-1)(54k^3-90k^2+33k-1)
$,

$ p_{\gamma_1,\beta_2}^{\gamma_2}=(2k-1)(3k-1)(18k^2-18k+5)
$,
\quad
$ p_{\gamma_1,\beta_1}^{\gamma_2}=(3k-2)(2k-1)(18k^2-18k+5)
$,

$ p_{\alpha_1,\gamma_2}^{\gamma_2}=k(54k^3-72k^2+15k+4)
$,
\qquad\qquad
$ p_{\alpha_2,\gamma_1}^{\gamma_2}=(18k^2-18k+5)(3k-1)(k-1)
$,

$ p_{\gamma_2,\beta_1}^{\gamma_2}=(6k-5)k(18k^2-18k+5)
$,
\qquad\qquad
$ p_{\alpha_1,\gamma_1}^{\gamma_2}=(3k-2)k(18k^2-18k+5)
$.\\

\noindent
In above $p_{a,b}^{c}
=p_{b,a}^{c}$ holds for any $a,b,c\in
\{\alpha_i,\beta_j, \gamma_k\mid
i,\ j=0,1,2, k=1,2\}$.\\

\section{Appendix III}
Intersection numbers of the coherent configuration
attached to the Euclidean $4$-design
supported by 2 concentric spheres satisfying 
$N_2\geq N_1\geq n+2$.

$A(X_1,X_1)=\{\alpha_1,\alpha_2\}$,
$A(X_2,X_2)=\{\beta_1,\beta_2\}$,
$A(X_1,X_2)=\{\gamma_1,\gamma_2\}$,

$\gamma_1\gamma_2=-\frac{1}{n}$,
$\alpha_2=-\frac{n\alpha_1-n+N_1}
{n((N_1-1)\alpha_1+1)}$,
$\beta_2=-\frac{n\beta_1-n+N_1}
{n((N_1-1)\beta_1+1)}$.\\

$p_{\alpha_1,\alpha_2}^{\alpha_0}
=p_{\alpha_2,\alpha_1}^{\alpha_0}
=p_{\gamma_2,\gamma_1}^{\alpha_0}
=p_{\gamma_1,\gamma_2}^{\alpha_0}
=0$,\quad
$p_{\alpha_2,\alpha_2}^{\alpha_0}=
N_1-p_{\alpha_1,\alpha_1}^{\alpha_0}-1$,\quad
$p_{\gamma_2,\gamma_2}^{\alpha_0}=
N_2-p_{\gamma_1,\gamma_1}^{\alpha_0}$,\\

$p_{\alpha_1,\alpha_1}^{\alpha_0}
= \frac{(N_1-1)\alpha_2+1}{(\alpha_2-\alpha_1)}$,
\quad
$p_{\gamma_1,\gamma_1}^{\alpha_0}
=\frac{N_2}{1+n\gamma_1^2}$.\\

$p_{\beta_1,\beta_2}^{\beta_0}
=p_{\beta_2,\beta_1}^{\beta_0}
=p_{\gamma_2,\gamma_1}^{\beta_0}
=p_{\gamma_1,\gamma_2}^{\beta_0}
=0$,\quad
$p_{\beta_2,\beta_2}^{\beta_0}=
N_2-p_{\beta_1,\beta_1}^{\beta_0}-1$,
\quad
$p_{\gamma_2,\gamma_2}^{\beta_0}=
N_1-p_{\gamma_1,\gamma_1}^{\beta_0}$,

$p_{\gamma_1, \gamma_1}^{\beta_0}
= \frac{N_1}{n\gamma_1^2+1}$,\quad
$p_{\beta_1, \beta_1}^{\beta_0}
= \frac{(N_2-1)\beta_2+1}{\beta_2-\beta_1}
$.\\
\begin{eqnarray}
&&p_{\alpha_1, \alpha_2}
^{\alpha_1} =
p_{\alpha_2, \alpha_1}^{\alpha_1} =
\frac{ n(1-\alpha_1)
(N_1\alpha_1-\alpha_1+1)^2
(n\alpha_1+N_1-n)}
{\bigg(N_1-n+2n\alpha_1+n(N_1-1)\alpha_1^2
 \bigg)^2},\nonumber\\
 &&p_{\alpha_1, \alpha_1}^{\alpha_1}=
 \nonumber\\
 &&\qquad \frac{\bigg(n(N_1-1)(N_1-2n-1)\alpha_1^3-3n^2\alpha_1^2-3n\alpha_1
 +(N_1-n-2)(N_1-n)\bigg)N_1}
 {\bigg(N_1-n+2n\alpha_1+n(N_1-1)\alpha_1^2
 \bigg)^2},
 \nonumber\\
  &&p_{\alpha_2, \alpha_2}^{\alpha_1} 
  = \frac{n\alpha_1N_1(N_1\alpha_1-\alpha_1+1)^2(n\alpha_1+1)}
{ \bigg(N_1-n+2n\alpha_1+n(N_1-1)\alpha_1^2
 \bigg)^2},
\nonumber\\
&&p_{\gamma_1,\gamma_1}^{\alpha_1}
= \frac{N_2(n\gamma_1^2\alpha_1+1)}{(\gamma_1^2n+1)^2},\ 
p_{\gamma_2, \gamma_2}^{\alpha_1} =
 \frac{(\gamma_1^2n+\alpha_1)N_2\gamma_1^2n}{(\gamma_1^2n+1)^2},\
 p_{\gamma_1, \gamma_2}^{\alpha_1}
  =p_{\gamma_2, \gamma_1}^{\alpha_1}
  = \frac{(1-\alpha_1)n\gamma_1^2N_2}{(\gamma_1^2n+1)^2}.
  \nonumber
\end{eqnarray}
\begin{eqnarray}
&&p_{\alpha_1, \alpha_2} 
^{\alpha_2}= 
p_{\alpha_2, \alpha_1} ^{\alpha_2}= 
\frac{N_1^2\alpha_1(n\alpha_1+1)(N_1-n-1)}
{\bigg(N_1-n+2n\alpha_1+n(N_1-1)\alpha_1^2
 \bigg)^2},
\nonumber\\
&&p_{\alpha_1, \alpha_1} ^{\alpha_2}= 
\frac{ N_1(1-\alpha_1)(N_1-n-1)(n\alpha_1-n+N_1)}{\bigg(N_1-n+2n\alpha_1+n(N_1-1)\alpha_1^2
 \bigg)^2},\nonumber\\
&&p_{\alpha_2, \alpha_2} ^{\alpha_2}=
\frac{(N_1-1)\alpha_1+1}
{\bigg(N_1-n+2n\alpha_1+n(N_1-1)\alpha_1^2
 \bigg)^2}\bigg(n^2(N_1^2-3N_1+2)\alpha_1^3\nonumber\\
&& 
+3n^2(N_1-2)\alpha_1^2
+3n(2n-N_1)\alpha_1-2n^2+3nN_1-N_1^2\bigg),
\nonumber\\
&&p_{\gamma_1, \gamma_2} ^{\alpha_2}= 
p_{\gamma_2, \gamma_1} ^{\alpha_2}= 
\frac{(n\alpha_1+1)\gamma_1^2N_2N_1}{(\gamma_1^2n+1)^2(N_1\alpha_1-\alpha_1+1)},
\nonumber\\
&&p_{\gamma_1, \gamma_1} ^{\alpha_2}=
\frac{  -N_2(-N_1\alpha_1+N_1\gamma_1^2+n\gamma_1^2\alpha_1+\alpha_1-\gamma_1^2n-1)}{(\gamma_1^2n+1)^2(N_1\alpha_1-\alpha_1+1)},
\nonumber\\
&&p_{\gamma_2, \gamma_2} ^{\alpha_2}= \frac{N_2\gamma_1^2(-N_1-n\alpha_1+n+n^2\alpha_1N_1\gamma_1^2-n^2\alpha_1\gamma_1^2
+\gamma_1^2n^2)}
{(\gamma_1^2n+1)^2(N_1\alpha_1-\alpha_1+1)}.\nonumber
\end{eqnarray}
\begin{eqnarray}
&&p_{\beta_1, \beta_2}^{\beta_1}=
p_{\beta_2, \beta_1}^{\beta_1}
= \frac{n(1-\beta_1)(1+N_2\beta_1-\beta_1)^2(-n+N_2+n\beta_1)}
{\bigg(N_2-n+2n\beta_1+n(N_2-1)\beta_1^2\bigg)^2}\nonumber\\
&&p_{\beta_1, \beta_1}^{\beta_1}=
\frac{N_2}{\bigg(N_2-n+2n\beta_1+n(N_2-1)\beta_1^2\bigg)^2}
\bigg(n(N_2-1)(N_2-2n-1)\beta_1^3\nonumber\\
&&
-3n^2\beta_1^2-3n\beta_1+(N_2-n-2)(N_2-n)\bigg)
\nonumber\\
&&p_{\beta_2, \beta_2}^{\beta_1}=
\frac{ N_2n\beta_1(1+N_2\beta_1-\beta_1)^2(n\beta_1+1)}{\bigg(N_2-n+2n\beta_1+n(N_2-1)\beta_1^2\bigg)^2},\
p_{\gamma_1, \gamma_2}^{\beta_1}=
p_{\gamma_2, \gamma_1}^{\beta_1}=
 \frac{(1-\beta_1)N_1\gamma_1^2n}{(\gamma_1^2n+1)^2},
 \nonumber\\
&&p_{\gamma_1, \gamma_1}^{\beta_1}=
\frac{ N_1(\gamma_1^2n\beta_1+1)}{(\gamma_1^2n+1)^2},\ 
p_{\gamma_2, \gamma_2}^{\beta_1}=
\frac{ N_1\gamma_1^2n(\beta_1+\gamma_1^2n)}{(\gamma_1^2n+1)^2}.\nonumber
\end{eqnarray}
\begin{eqnarray}
&&
p_{\beta_1, \beta_2}^{\beta_2} = 
p_{\beta_2, \beta_1}^{\beta_2}=
\frac{N_2^2\beta_1(n\beta_1+1)(N_2-n-1)}
{\bigg(N_2-n+2n\beta_1+n(N_2-1)\beta_1^2\bigg)^2},
\nonumber\\
&&p_{\beta_1, \beta_1}^{\beta_2} 
= \frac{(-1+\beta_1)(n-N_2+1)(-n+N_2+n\beta_1)N_2}{\bigg(N_2-n+2n\beta_1+n(N_2-1)\beta_1^2\bigg)^2},
\nonumber\\
&&p_{\beta_2, \beta_2}^{\beta_2} 
=\frac{(1+N_2\beta_1-\beta_1)}
{\bigg(N_2-n+2n\beta_1+n(N_2-1)\beta_1^2\bigg)^2} 
\bigg(n^2(N_2^2-3N_2+2)\beta_1^3
+3n^2\beta_1^2(N_2-2)\nonumber\\
&&\qquad-3n(N_2-2n)\beta_1-2n^2+3N_2n-N_2^2\bigg),\nonumber\\
&&p_{\gamma_1, \gamma_2}^{\beta_2} =
p_{\gamma_2, \gamma_1}^{\beta_2}=
\frac{ (n\beta_1+1)\gamma_1^2N_1N_2}{(\gamma_1^2n+1)^2(1+N_2\beta_1-\beta_1)},
\nonumber\\
&&
p_{\gamma_1, \gamma_1}^{\beta_2} = 
\frac{\big((n-N_2-n\beta_1)\gamma_1^2
+N_2\beta_1+1-\beta_1\big)N_1}{(\gamma_1^2n+1)^2(1+N_2\beta_1-\beta_1)},
\nonumber\\
&&p_{\gamma_2, \gamma_2}^{\beta_2} 
=\frac{ ((n^2\beta_1N_2-n^2\beta_1+n^2)\gamma_1^2
-n\beta_1+n-N_2)N_1\gamma_1^2}{(\gamma_1^2n+1)^2(1+N_2\beta_1-\beta_1)}.\nonumber
\end{eqnarray}
\begin{eqnarray}
&&p_{\gamma_2, \beta_2}
^{\gamma_1} =
\frac{ N_2\gamma_1^2n(1+N_2\beta_1-\beta_1)(n\beta_1+1)}
{\bigg(N_2-n+2n\beta_1+n(N_2-1)\beta_1^2\bigg)(\gamma_1^2n+1)},
\nonumber\\
&&p_{\alpha_2, \gamma_2}^{\gamma_1}
= \frac{N_1\gamma_1^2n(N_1\alpha_1-\alpha_1+1)(n\alpha_1+1)}
{\bigg(N_1-n+2n\alpha_1+n(N_1-1)\alpha_1^2\bigg)(\gamma_1^2n+1)},
\nonumber\\
&&p_{\gamma_1, \beta_1}^{\gamma_1} = 
\frac{(\gamma_1^2n\beta_1+1)(N_2-n-1)N_2}
{\bigg(N_2-n+2n\beta_1+n(N_2-1)\beta_1^2\bigg)(\gamma_1^2n+1)},
\nonumber\\
&&
p_{\gamma_2, \beta_1}^{\gamma_1}
 =\frac{ N_2\gamma_1^2n(1-\beta_1)(N_2-n-1)}
 {\bigg(N_2-n+2n\beta_1+n(N_2-1)\beta_1^2\bigg)(\gamma_1^2n+1)},
 \nonumber\\
&&p_{\alpha_1, \gamma_1}^{\gamma_1} =
\frac{ N_1(n\gamma_1^2\alpha_1+1)(N_1-n-1)}
{\bigg(N_1-n+2n\alpha_1+n(N_1-1)\alpha_1^2\bigg)(\gamma_1^2n+1)},
\nonumber\\
&&p_{\alpha_1, \gamma_2}^{\gamma_1} 
= \frac{N_1\gamma_1^2n(1-\alpha_1)(N_1-n-1)}
{\bigg(N_1-n+2n\alpha_1+n(N_1-1)\alpha_1^2\bigg)(\gamma_1^2n+1)},
\nonumber\\
&&p_{\gamma_1, \beta_2}^{\gamma_1} 
=\frac{
 n(1+N_2\beta_1-\beta_1)(
 (n-n\beta_1-N_2)\gamma_1^2+N_2\beta_1+1-\beta_1)}{\bigg(N_2-n+2n\beta_1+n(N_2-1)\beta_1^2\bigg)(\gamma_1^2n+1)},
 \nonumber\\
&&p_{\alpha_2, \gamma_1}^{\gamma_1} = 
\frac{n(N_1\alpha_1-\alpha_1+1)
(N_1\alpha_1-\alpha_1+1
+(n-N_1-n\alpha_1)\gamma_1^2)}
{\bigg(N_1-n+2n\alpha_1+n(N_1-1)\alpha_1^2\bigg)(\gamma_1^2n+1)}.\nonumber
\end{eqnarray}
\begin{eqnarray}
&&p_{\gamma_2, \beta_2}^{\gamma_2} 
=\frac{ (1+N_2\beta_1-\beta_1)
((n^2\beta_1N_2-n^2\beta_1+n^2)\gamma_1^2
-n\beta_1+n-N_2)}
{\bigg(N_2-n+2n\beta_1+n(N_2-1)\beta_1^2\bigg)(\gamma_1^2n+1)},
\nonumber\\
&&p_{\alpha_2, \gamma_2}^{\gamma_2}
 =\frac{ (N_1\alpha_1-\alpha_1+1)
 (-N_1-n\alpha_1+n
 +(n^2\alpha_1N_1-n^2\alpha_1+n^2)\gamma_1^2)}{\bigg(N_1-n+2n\alpha_1+n(N_1-1)\alpha_1^2\bigg)(\gamma_1^2n+1)},\nonumber\\
&&p_{\gamma_1, \beta_2}^{\gamma_2} 
=\frac{ (1+N_2\beta_1-\beta_1)(n\beta_1+1)N_2}{\bigg(N_2-n+2n\beta_1+n(N_2-1)\beta_1^2\bigg)(\gamma_1^2n+1)},
\nonumber\\
&&p_{\gamma_1, \beta_1}^{\gamma_2} =
\frac{ N_2(-1+\beta_1)(n-N_2+1)}
{\bigg(N_2-n+2n\beta_1+n(N_2-1)\beta_1^2\bigg)(\gamma_1^2n+1)},\nonumber\\
&&p_{\alpha_1, \gamma_2}^{\gamma_2} =
\frac{ (\gamma_1^2n+\alpha_1)(N_1-n-1)N_1}
{\bigg(N_1-n+2n\alpha_1+n(N_1-1)\alpha_1^2\bigg)(\gamma_1^2n+1)},
\nonumber\\
&&p_{\alpha_2, \gamma_1}^{\gamma_2} 
=\frac{ N_1(N_1\alpha_1-\alpha_1+1)(n\alpha_1+1)}
{\bigg(N_1-n+2n\alpha_1+n(N_1-1)\alpha_1^2\bigg)(\gamma_1^2n+1)},\nonumber\\
&&p_{\gamma_2, \beta_1}^{\gamma_2} = 
\frac{N_2(N_2-n-1)(\beta_1+\gamma_1^2n)}
{\bigg(N_2-n+2n\beta_1+n(N_2-1)\beta_1^2\bigg)(\gamma_1^2n+1)},\nonumber\\
&&p_{\alpha_1, \gamma_1}^{\gamma_2} = 
\frac{(1-\alpha_1)(N_1-n-1)N_1}
{\bigg(N_1-n+2n\alpha_1+n(N_1-1)\alpha_1^2\bigg)(\gamma_1^2n+1)}.
\nonumber
\end{eqnarray}

Acknowledgment:
The authors thank Masanobu Kaneko of Kyushu University for providing a proof of some results on
diophantine equations, i.e. Proposition \ref{pro:Kaneko1}
and Proposition \ref{pro:Kaneko2}.

\newpage
\noindent
2000 Mathematics
Subject Classification.

Primary: 05E99, Secondary:
05B99, 51M99, 62K99.\\

\noindent
Key words and Phrases.

Euclidean design, spherical design, 
association scheme,
coherent configuration,
cubature formula

\end{document}